\documentclass[11pt]{article}
\usepackage{amsmath,dsfont,amssymb,latexsym,bm,amsthm,amsfonts,float,graphicx,xcolor}

\theoremstyle{thmstyleone}
\newtheorem{theorem}{Theorem}
\newtheorem{lemma}{Lemma}
\newtheorem{remark}{Remark}

\newcommand{\norm}[1]{\left\Vert#1\right\Vert}
\newcommand{\abs}[1]{\left\vert#1\right\vert}
\def\E{\mathbb{E}}
\def\Cov{{\rm{Cov}}}
\def\Var{{\rm{Var}}}

\def\P{\mathbb{P}}
\def\B{\mathbb{B}}

\def\R{\mathds{R}}
\def\Z{\mathds{Z}}
\def\NN{\mathcal{N}}
\def\OO{\mathcal{O}}\def\QQ{\mathcal{Q}}

\def\UU{\mathcal{U}}
\def\JJ{\mathcal{J}}

\def\LL{\mathcal{L}}
\def\CC{\mathcal{C}}

\def\FF{\mathcal{F}}

\def\KK{\mathcal{K}}

\def\AA{\mathcal{A}}
\def\RR{\mathcal{R}}
\everymath{\displaystyle}

\theoremstyle{thmstyletwo}

\theoremstyle{thmstylethree}

\begin{document}

\begin{center}
{\LARGE\bf Two approaches of estimation for a tempered fractional Brownian motion with drift}
\\[1ex]
by
\\[1ex]
Mohamed El Omari$^1$ and Saralees Nadarajah$^{2,  *}$
\\
$^1$ESIM, Polydisciplinary Faculty of Sidi Bennour, Chouaib Doukkali University, El Jadida, 24000, Morocco, {\sf elomari.m@ucd.ac.ma}
\\
$^{2}$Department of Mathematics, University of Manchester, Manchester M13 9PL, UK, {\sf mbbsssn2@manchester.ac.uk}
\\
$^*$corresponding author
\end{center}
\vspace{1.5cm}
{\bf Abstract.}~~We investigate the problem of estimation for a tempered
fractional Brownian motion with drift of the form: $X_t=\mu t + \sigma B_{H,\lambda}(t)$, for all $t\geq 0$.
Based on discrete observations of $X$ and by using the method of moments with the ergodic theorem,
we construct estimators of the drift parameter $\mu\in\R$ and noise parameters  $ \sigma^2,H,\lambda \in(0,\infty)$.
We establish the asymptotic properties of the proposed estimators, including
convergence in probability, almost sure convergence, and asymptotic normality. Finally,
by using Gladyshev's results \cite{gladyshev1961new} we provide some consistent alternative estimators of   $\sigma^2$, $H$ and $\lambda$, when $\mu=0$ and $H\in (0,1)$.



\section{Introduction}
\label{Sec.1}

Tempered fractional Brownian motions (tfBms) of first or second kind
have been proposed by \cite{Meerschaert2013,sabzikar2018tempered} as an extension of the fractional Brownian motion (fBm).
They are constructed by exponentially tempering the power law kernel in the Mandelbrot and Van Ness representation of the fBm.
A tfBm $B_{H,\lambda}$ shares the Gaussianity and stationarity of increments with the fBm,
but is characterized by a mean squared displacement which decays rapidly to a constant,  
i.e., ${\rm{Var}} \left(B_{H,\lambda}(t)\right)=\CC +\OO  ( t^{H-1/2}e^{-\lambda t} )$
, as $t\rightarrow \infty$, for some constant $\CC >0$. The autocorrelation of its increments exhibits the semi-long range dependence,
that is, it behaves like a power over moderate scales, but quasi-exponentially over large scales.
TfBms are widely used to model wind speed data
summarized by the so-called Davenport spectrum (e.g., \cite{davenport1961spectrum}).
The latter is a modification of the classical Kolmogorov spectrum for the inertial range of turbulence at low frequencies.
TfBms have also been used to model geophysical flow data \cite{boniece2018tempered} and in finance \cite{zhang2017arbitrage}.

The theoretical aspects of tfBm have attracted considerable research interest (see \cite{azmoodeh2022does,liu2020numerical,Meerschaert2013,meerschaert2014stochastic,shen2017strong}
and references therein).
In contrast, literature on their statistical aspect remains limited to a few papers
\cite{boniece2021tempered,boniece2018tempered}, in which wavelet estimation procedures were performed.
Recently, least-squares estimators for the unknown drift parameters within Vasicek models driven by tfBms have been proposed
and examined in \cite{mishura2024asymptotic,mishura2025asymptotic}.
This work addresses the parameter estimation of a tfBm with drift of the form.
\begin{equation}
\label{Model}
\displaystyle
X_t = \mu t + \sigma B_{H,\lambda}(t)
\end{equation}
for all $t \geq 0$, where $B_{H,\lambda}$ is a tfBm as defined below.
Using the ergodic theorem and the method of moments (e.g., \cite{barboza2017parameter})
we derive suitable estimators for model \eqref{Model}.
We investigate both the weak and almost sure behaviour of the proposed estimators.
Consequently, our results provide alternatives to those given in \cite{boniece2021tempered,boniece2018tempered}.
When $H\in(0,1)$, we apply Gladyshev's results to construct some alternative estimators 
of the noise parameters $H$, $\sigma^2$, $\lambda$ and prove their consistency.
The evaluation of their limiting distributions requires more advanced tools and is beyond the scope of this work.

The rest of the paper is organized as follows.
In Section \ref{Sec.2} we introduce some basic and preliminary results.
Then we discuss the procedures of estimation to be used.
In Section \ref{Sec.3} we construct $\widehat{\mu}$ the estimator of the drift $\mu\in\R$ and investigate
its asymptotic behaviour (independently of the other parameters).
Section \ref{Sec.4} is devoted to the joint estimation of the noise
parameters  $(\sigma^2,H,\lambda)\in(0,\infty)^3$ by the method of moments (e.g., \cite{barboza2017parameter}).
In Section \ref{Sec.5} we use Gladyshev's results \cite{gladyshev1961new} to propose a multistage
estimation for the noise parameters $H$, $\sigma^2$, $\lambda$ under the condition $0<H<1$.
In this section, the consistency of the estimators of $H$ and $\sigma^2$ holds almost surely,
while the estimator of $\lambda$ is consistent in probability.
Section \ref{con} concludes the paper.
Some auxiliary and technical lemmas on the modified Bessel function $K_H(\cdot)$ with
the Birkhoff's ergodic theorem are gathered in the Appendix.

In what follows we work on a probability space $(\Omega,\FF,\P)$ and each process
should be understood as being defined with respect to its own natural filtration.
We are concerned with four modes of convergence: almost surely $\overset{\mathbb{P}-a.s}{\longrightarrow}$,
in the mean square sense $\overset{L^2(\Omega)}{\longrightarrow}$,
in probability $\overset{\mathbb{P}}{\longrightarrow}$ and in law $\overset{\mathcal{L}}{\longrightarrow}$.
The following symbols  $\OO(\cdot)$ and $o(\cdot)$ denote the usual  \lq\lq big-oh\rq\rq and  \lq\lq small-oh\rq\rq notations, respectively.

\section{Preliminary results}
\label{Sec.2}

Let $\left(B(t):t\in\R\right)$ be a two-sided Brownian motion,
that is, a process with stationary independent increments
such that $B(t)$ has a Gaussian distribution with mean zero and variance $\abs{t}$ for all $t\in\R$.
Define an independently scattered Gaussian random
measure $B(dx)$ with control measure $m (dx)=dx$ by setting $B([a,b])=B(b)-B(a)$ for any $a,b\in\R$
with $a<b$ and then extending to all Borel sets with finite Lebesgue measure.
It is well known that the Wiener integral $\int f(x) B(dx)$ exists whenever $\int f(x)^2 dx <\infty$.
Many famous fractional processes were constructed by choosing a suitable integrand $f(\cdot)$.
For instance, the fractional Brownian motion $B_H$ was introduced by setting $f(x) = C_H \left( (t-x)_+^{H-1/2}-(-x)_+^{H-1/2} \right)$,
where $C_H>0$ is some  normalizing constant.
In \cite{Meerschaert2013}, the integrand in the form of $B_H$ is exponentially tempered to get the tfBm 
\begin{equation*}
\label{tfBm}
\displaystyle
B_{H,\lambda}(t) = \int \left[e^{-\lambda(t-x)_+}(t-x)_+^{H-1/2}-e^{-\lambda(-x)_+}(-x)_+^{H-1/2} \right] B(dx)
\end{equation*}
for all $t \in R$.
$B_{H,\lambda}$ exhibits many interesting properties
such as stationarity of increments and semi-long range dependence.
It is worth mentioning
that this process is not self-similar, but instead the following
scaling property: $B_{H,\lambda}(ct)\triangleq c^H B_{H,c\lambda}(t)$ holds for any $c>0$.
This property implies
\begin{eqnarray}
\label{VartfBm}
\displaystyle
\E \left[ B_{H,\lambda}(\vert t \vert)^2 \right]
&=&
\displaystyle
t^{2H} \E \left[ B_{H,  \lambda \vert t \vert}(1)^2 \right]=: t^{2H} C_t^2,
\\
\nonumber
&&
\displaystyle
\longrightarrow 2\Gamma(2H)/(2\lambda)^{2H}
\end{eqnarray}
as $t \to \infty$, where 
\begin{eqnarray*}
\displaystyle
C_t^2  =  \frac {2\Gamma (2H)}{\left(2\lambda \vert t \vert \right)^{2H}} -
\frac {2\Gamma \left(H+1/2\right)}{\sqrt{\pi} \left(2\lambda \vert t \vert \right)^{H}}K_H \left(\lambda\vert t\vert\right)~~\mbox{ and }C_0^2=0,
\end{eqnarray*}
where $K_H$ denotes the modified Bessel function of second kind.
The covariance of $B_{H,\lambda}$ can be expressed as
\begin{equation}\label{COVFUN}
\displaystyle
\Cov \left(B_{H,\lambda}(t),B_{H,\lambda}(s)\right) = \frac {1}{2}\left(\abs{t}^{2H}C^2_t+\abs{s}^{2H}C^2_s-\abs{t-s}^{2H}C^2_{t-s}\right)
\end{equation}
for all $t, s \geq 0$.
When $H\in (0,1)$, the tfBm $B_{H,\lambda}$ becomes a quasi-hélix process in
the sense of  \cite{kahane1980helices}, and the H\"{o}lder continuity follows with $H$ as critical exponent.
For $H>1$, the sample paths of $B_{H,\lambda}$ are  continuously differentiable.
For more details on properties of the tfBm we refer to \cite{Meerschaert2013,azmoodeh2022does}.

Given a function $g:\R\longrightarrow \R$, we define its $k$th order increments with step $h\in\R$ as 
\begin{equation*}
\label{increments}
\displaystyle
\Delta^{(k)}_h g(t) := \sum_{j=0}^k (-1)^{k-j}
\begin{pmatrix}
k
\\
j
\end{pmatrix} g(t+hj)
\end{equation*}
for all $t\in\R$ and  $~k\in\Z_+$ with the convention $\Delta^{(0)}_h g(t)=g(t)$.
We say that a process $Y=\left(Y(t):t\in\R\right)$ is $k$-stationary if for every $h\in\R$ the processes
$\left(\Delta^{(r)}_h Y(t):t\in\R\right)$, $r=0,\ldots,(k-1)$ are non-stationary, but  $\left(\Delta^{(k)}_h Y(t):t\in\R\right)$ exhibits stationarity.
The process $X$ defined in \eqref{Model} is $1$-stationary and stationarity holds also for any increments of higher order $k \geq 2$.
Consider the following statistics to be used later
\begin{equation}
\label{eq3}
\displaystyle
\overline{X}_N^{p,q}(h) :=\frac {1}{N} \sum_{k=0}^{N-1} \left(\Delta^{(p)}_h X_t\mid_{t=kh} \right)^q,~~N=1,2,\ldots,
\end{equation}
for fixed $h>0$ and $p,q\in\Z_+^*$.
Then the following results are obtained.

\begin{lemma}
\label{LEM1}
We have  
\begin{align*}
\displaystyle
\overline{X}_N^{p,q}(h)\overset{\P-a.s}{\longrightarrow} \E \left(\Delta^{(p)}_h X_t\mid_{t=0} \right)^q =:\overline{X}_\infty^{p,q}(h)
\end{align*}
as $N \to \infty$.
In particular,
\begin{equation}
\label{eq4}
\begin{cases}
\displaystyle
\overline{X}_\infty^{1,1}(h)=\mu h,~~\overline{X}_\infty^{1,2}(h)=\mu^2h^2+\sigma^2  h^{2H}C_h^2,
\\
\displaystyle
\overline{X}_\infty^{2,2}(h)=\sigma^2 h^{2H} \left( 4C_h^2-2^{2H}C_{2h}^2\right),
\end{cases}
\end{equation}
\end{lemma}

\begin{proof}
\label{proof-LEM1}
Let $Z_h(t):=\Delta^{(1)}_h X_t$.
Because of stationarity of $Z_h$, we have
\begin{eqnarray}
\nonumber
&&
\displaystyle
{\rm{Cov}}\left(Z_h(s+\tau ),Z_h(s)\right) = {\rm{Cov}}\left(Z_h(\tau),Z_h(0)\right)
\\
\nonumber
&&=
\displaystyle
\sigma^2 \E \left\lbrace B_{H,\lambda}(h)\left(B_{H,\lambda}(\tau +h)-B_{H,\lambda}(\tau)\right) \right\rbrace
\\
\label{eq4.9}
&&=
\displaystyle
\frac{\sigma^2}{2}
 \left( \abs{\tau+h}^{2H}C^2_{\tau+h}+\abs{\tau-h}^{2H}C^2_{\tau-h}-2\abs{\tau}^{2H}C^2_\tau \right),
\\
\label{eq5}
&&=
\displaystyle
\frac {\sigma^2\Gamma \left(H+1/2\right)}{\sqrt{\pi} 2^H\lambda^{2H}}
\left[ 2 g_H \left(\lambda \tau  \right) -
g_H \left(\lambda (\tau +h) \right) - g_H \left(\lambda ( \tau -h) \right) \right]
\\
\nonumber
&&
\displaystyle
\longrightarrow 0
\end{eqnarray}
as $\tau \to \infty$, where $g_H(z)=\vert z \vert^H K_H(\vert z \vert)$.
The above convergence is justified by the fact that $g_H(z)\sim \sqrt{\frac {\pi}{2}}z^{H-1/2}e^{-z}$, as $z\rightarrow\infty$ \cite[(A.5)]{gaunt2014inequalities}.
 For each $p\in\Z_+^*$  the covariance of $\Delta^{(p)}_h X$ is continuous and satisfies
\begin{eqnarray}
\label{eq6}
&&
\displaystyle
{\rm{Cov}}\left(\Delta^{(p)}_h X(s+\tau),
\Delta^{(p)}_h X(s)\right) =
\sum_{j,r=0}^{p-1} d_{j,r}  {\rm{Cov}} \left(Z_h \left(\tau +h(j-r)\right), Z_h(0)\right)
\\
\nonumber
&&
\displaystyle
\longrightarrow 0
\end{eqnarray}
as $\tau \to \infty$, where $d_{j,r} = (-1)^{j+r} \begin{pmatrix}p-1\\j\end{pmatrix} \begin{pmatrix}p-1\\r\end{pmatrix} $.
Combining the previous asymptotics with Gaussianity and stationarity
of the sequence $\lbrace \Delta^{(p)}_h X_t\mid_{t=kh}:k=0,1,2,\ldots \rbrace $ yields its ergodicity 
(e.g., \cite[Theorem 5, p. 144]{gikhman2004theory}).
Hence, by Lemma \ref{LEM1-App} we can justify the first part of Lemma \ref{LEM1}.
The proof of \eqref{eq4} follows by elementary calculus with the use of \eqref{VartfBm}.  
\end{proof}

Lemma \ref{LEM1} allows us to construct a vector of
statistics \eqref{eq3} (say $J_N\in\R^4$) that converges towards a function of the parameters (say $F(\mu,\sigma,H,\lambda)\in\R^4$).
Hence, an appropriate estimator of $(\mu,\sigma,H,\lambda)$ can be defined as $F^{-1}\left(J_N\right)$.
This approach has been successfully applied to various mixed models (see \cite{ralchenko2023parameter,ralchenko2023asymptotic,el2024statistical}).
However, finding such an inverse is often difficult, if not impossible.
To overcome this problem we estimate the parameters in two steps:
First, we estimate the drift parameter $\mu$ by $\overline{X}_N^{1,1}(h)/h$. This estimator can be seen as a discrete form of $\widehat{\mu}_t:=X_t/t$.
Second, we follow \cite{barboza2017parameter} to perform the joint estimation of $\left(\sigma^2,H,\lambda\right)$
by the method of generalized moments.
Note that a joint estimation of all the parameters is not
possible by this method, as the identifiability condition (see \cite{barboza2017parameter}) does not hold.
Finally, we apply Gladyshev's results on  limit theorems for processes
with Gaussian increments to $X$ and derive consistent (alternative)
estimators of the noise parameters denoted in the following order
as $\widetilde{H}_{n}$, $\widetilde{\sigma}^2_{m,n}$ and $\widetilde{\lambda}_{m,n,h}$.
This order is crucial for the implementation issue.
In fact, we can write the latter two estimators as $\widetilde{\sigma}^2_{m,n}=F_m (\widetilde{H}_{n},X )$ and
$\widetilde{\lambda}_{m,n,h}=G_h ( \widetilde{H}_{n},\widetilde{\sigma}^2_{m,n},X )$
for some continuous functions $F_m$ and $G_h$.

\section{Estimation of the drift parameter $\mu$}
\label{Sec.3}

\begin{theorem}
\label{TH1}
For a fixed $h>0$, the following statistic $\widehat{\mu}:=\overline{X}_N^{1,1}(h)/h$
converges towards $\mu$, almost surely and in $L^2(\Omega)$, as $N\rightarrow\infty$.
Furthermore, $\widehat{\mu}$ is unbiased and satisfies the asymptotic normality i.e.,
\begin{equation}
\label{eq7}
\displaystyle
N\left(\widehat{\mu}-\mu\right)\overset{\LL}{\longrightarrow} \NN \left(0, \frac{2\sigma^2\Gamma (2H)}{h^2(2\lambda )^{2H}}\right),~~\mbox{as }N \to \infty
\end{equation}  
\end{theorem}

\begin{proof}
\label{proof-TH1}
The unbiasedness of $\widehat{\mu}$ is trivial and its strong
consistency follows immediately from Lemma \ref{LEM1}.
Let us evaluate the variance of $\widehat{\mu}$.
Consider $Z_h(\cdot)$ as defined in the proof of Lemma \ref{LEM1} and set $\phi_h(\tau):=\abs{h\tau}^{2H}C^2_{h\tau}$.
In view of \eqref{eq4.9}, we have
\begin{eqnarray*}
&&
\displaystyle
N^2 h^{2}\Var  (\widehat{\mu} ) =  \Var   ( \sum_{k=0}^{N-1} Z_h(hk) )  =
\Cov  \left( \sum_{k=0}^{N-1} Z_h(hk), \sum_{j=0}^{N-1} Z_h(hj) \right)
\\
&&=  
\E   \left[ \sum_{k=0}^{N-1} Z_h(hk)\sum_{j=0}^{N-1} Z_h(hj) \right] - \left(\E \sum_{k=0}^{N-1} Z_h(hk) \right) \left(\E\sum_{j=0}^{N-1} Z_h(hj) \right)
\\
&&=
\displaystyle
\sum_{k,j=0}^{N-1} \Cov \left(  Z_h \left(h(k-j)\right), Z_h(0)\right),~~~~(\mbox{By stationarity of }Z_h),
\\
&&=
\displaystyle
\frac{\sigma^2}{2}
 \sum_{k,j=0}^{N-1} \left[ \phi_h(  k-j+1)+\phi_h (  k-j-1)-2\phi_h(  k-j) \right],~~~~(\mbox{By \eqref{eq4.9} }),
\\
&&
=\frac{\sigma^2}{2}\left\lbrace \sum_{k=j} \left[ \phi_h(  k-j+1)+\phi_h (  k-j-1)-2\phi_h(  k-j) \right] \right.
 \\
&&   \left.~~~~~~~~~~~~ + \sum_{  k\neq j  }  \left[ \phi_h(  k-j+1)+\phi_h (  k-j-1)-2\phi_h(  k-j) \right] \right\rbrace \\
 &&=
\displaystyle
 \sigma^2\left\lbrace \sum_{k=0}^{N-1}  \left(\phi_h(1)-\phi_h(0)\right)  + \sum_{  j<k } \left( \phi_h(  k-j+1)+\phi_h (  k-j-1)-2\phi_h(  k-j) \right)\right\rbrace
\\
&&=
\displaystyle
 \sigma^2 \left(  N \phi_h(1) + \sum_{k=1}^{N-1}\sum_{r=1}^{k}  \left( \phi_h(  r+1)-\phi_h( r) \right)+
\sum_{k=1}^{N-1}\sum_{r=1}^{k} \left( \phi_h ( r-1)-\phi_h( r) \right) \right)
\\
&&=
\displaystyle
 \sigma^2 \left(  N \phi_h(1) + \sum_{k=1}^{N-1}  \left( \phi_h(  k+1)-\phi_h(1) \right)+   \sum_{k=1}^{N-1} \left(\phi_h (0)-\phi_h(k) \right) \right)
\\
&&=
\displaystyle
 \sigma^2 \left(  N \phi_h(1) -(N-1) \phi_h(1) +\sum_{k=1}^{N-1}   (\phi_h(  k+1)-\phi_h(k) ) \right)
\\
&&=
\displaystyle
   \sigma^2\abs{hN}^{2H}C_{hN}^2 =\sigma^2\left(\frac {2\Gamma (2H)}{\left(2\lambda \right)^{2H}} -
\frac {2\Gamma \left(H+1/2\right)}{\sqrt{\pi} 2^H\lambda^{2H}}  g_H \left(\lambda hN\right) \right)
\end{eqnarray*}
At this stage we can assert that $\Var \left(hN \left(\widehat{\mu}-\mu\right)\right)\sim     2\sigma^2\Gamma (2H)/(2\lambda )^{2H}$, as $N\rightarrow\infty$.
This implies the consistency of $\widehat{\mu}$ in $L^2(\Omega)$.
Now, due to Gaussianity of $\widehat{\mu}$, equation \eqref{eq7} holds, and the proof is complete.  
\end{proof}

\section{Estimation of the noise parameters $(\sigma^2,H,\lambda )$}
\label{Sec.4}

Consider the process $Y_h=\lbrace \Delta_h^{(2)} X_t:t\geq 0\rbrace$ with $h>0$ to be specified later.
Clearly, this process is stationary Gaussian with zero mean.
It admits the following harmonizable representation (e.g., \cite[(3.1)]{Meerschaert2013})
\begin{eqnarray*}
\displaystyle
Y_h(t)
&=&
\displaystyle
\sigma \left( B_{H,\lambda}(t+2h)-2 B_{H,\lambda}(t+h)+B_{H,\lambda}(t)\right)
\\
&=&
\displaystyle
\frac {\sigma \Gamma (H+1/2)}{\sqrt{2\pi}}\int_\R \frac { \left(e^{-i2\omega h}-2e^{-i\omega h}+1\right)}{(\lambda-i\omega)^{H+1/2}}e^{-i\omega t} \widetilde{B}(d\omega),~~~(i=\sqrt{-1}),
\end{eqnarray*}
where $\widetilde{B}(d\omega)$ is a complex-valued Gaussian random measure satisfying
$\overline{\widetilde{B}(d\omega)} =  \widetilde{B}(-d\omega)$ and $\E\vert \widetilde{B}(d\omega)\vert^2 =  d\omega$.
We derive the spectral density function of $Y_h$ \cite[§. 9]{yaglom1987correlation} as
\begin{equation}
\label{sdf}
\displaystyle
f_\theta(\omega) = \frac {4\sigma^2 \Gamma (H+1/2)^2 \left(1-\cos (\omega h)\right)^2}{2\pi \left(\lambda^2+\omega^2\right)^{H+1/2}},~~\theta=\left( \sigma^2,H,\lambda \right).
\end{equation}
Its autocovariance structure is described by
\begin{eqnarray}
\nonumber
&&
\displaystyle
\gamma_\theta (\tau) := \Cov \left(Y_h(\tau),Y_h(0)\right)
\\
\label{eq7.9}
&&=
\displaystyle 
\frac{\sigma^2}{2}
\left\lbrace -\phi(\tau+2h)- \phi(\tau-2h)+4 \phi(\tau+h)+4 \phi(\tau-h) -6 \phi(\tau)\right\rbrace
\\ 
\label{eq8} &&=
\displaystyle
c \left\lbrace g_H \left(\lambda(\tau+2h)\right) - 4g_H\left(\lambda(\tau+h)\right)+6g_H(\lambda\tau)-
4g_H\left(\lambda(\tau-h)\right) +g_H\left(\lambda(\tau-2h)\right)\right\rbrace,
\end{eqnarray}
where $\phi(\tau)=\abs{\tau}^{2H} C^2_\tau$ and 
\begin{equation}
\label{eq9}
\displaystyle 
c= \frac {\sigma^2\Gamma \left(H+1/2\right)}{\sqrt{\pi} 2^H\lambda^{2H}}=:\sigma^2\kappa ~\mbox{ and }~g_H(z)=\vert z\vert^H K_H\left(\vert z\vert\right).
\end{equation}
Both \eqref{eq7.9} and \eqref{eq8}-\eqref{eq9} are of great importance.
The latter is valid if all the arguments used in $g_H(\cdot)$ are strictly positive.    For \eqref{eq8} to be valid for any argument including the origin we shall add the convention: $g_H(0)=2^{H-1}\Gamma(H)$. Now, Instead of \eqref{COVFUN} the covariance  can be rewritten  as
\begin{equation*} 
\Cov \left(B_{H,\lambda}(t),B_{H,\lambda}(s)\right) = \frac {1}{2}\left\lbrace c_1 +c_2\left[ g_H \left(\lambda(t-s)\right)-g_H\left(\lambda t\right)-g_H\left(\lambda s\right)\right]\right\rbrace,~~\forall t,s\in\R,
\end{equation*}where 
\begin{eqnarray*}
&& c_1 =\frac {2\Gamma (2H)}{\left(2\lambda \right)^{2H}}  
~~\mbox{and}~~c_2 =\frac {2\Gamma \left(H+1/2\right)}{\sqrt{\pi}  2^H\lambda^{2H}}.
\end{eqnarray*}In particular, 
$\E\abs{B_{H,\lambda}(t)}^2=c_1-c_2g_H(\lambda t)$, 
$\forall t\in\R$. Note that  $t=s=0$ implies
 $\E\abs{B_{H,\lambda}(0)}^2=c_1-c_2g_H(0)=0$,
  which shows that $B_{H,\lambda}$ starts from the origin. Note that $\phi(0)=\abs{0}^{2H}C_0^2=0$ but $g_H(0)\neq 0$. Let $\alpha>0$ (to be specified later) and set $\varrho_{\theta,3}(\alpha):=\left(\gamma_\theta (0),\gamma_\theta (\alpha),\gamma_\theta (2\alpha)\right)'$.
As the step $h$ is arbitrary we simplify our computations by setting $h=\alpha$ to obtain
\begin{equation*}
\displaystyle
\varrho_{\theta,3}(\alpha) =
c\begin{pmatrix} 
\displaystyle
2  \left\lbrace g_H(2\lambda\alpha)-4g_H(\lambda\alpha)\right\rbrace +6g_H(0)
\\
\displaystyle
g_H(3\lambda\alpha)-4g_H(2\lambda\alpha)+7g_H(\lambda\alpha)-4g_H(0)
\\
\displaystyle
g_H(4\lambda\alpha)-4g_H(3\lambda\alpha)+6g_H(2\lambda\alpha)-4g_H(\lambda\alpha) +g_H(0)
\end{pmatrix},~~~~\theta=\left( \sigma^2,H,\lambda \right)~,
\end{equation*}
with the Jacobian matrix with respect to $\theta$ given as
\begin{equation}
\label{eq10}
\displaystyle
\nabla \varrho_{\theta,3}(\alpha) = \begin{pmatrix}
\displaystyle
\Sigma_{1,1}~~ & \Sigma_{1,2}~~ & \Sigma_{1,3}
\\
\displaystyle
\Sigma_{2,1}~~ & \Sigma_{2,2}~~ & \Sigma_{2,3}
\\
\displaystyle
\Sigma_{3,1}~~ & \Sigma_{3,2}~~ & \Sigma_{3,3} 
\end{pmatrix},
\end{equation}
and
\begin{eqnarray}
\nonumber &&
\Sigma_{1,1} = 2\kappa(y-4x+3u_0)\\ 
&&
\label{pre-xyzt15} \Sigma_{2,1} =\kappa(z-4y+7x-4u_0)\\
\nonumber &&
\Sigma_{3,1} =\kappa(t-4z+6y-4x+u_0)\\
&&
\nonumber\\
\nonumber &&\Sigma_{1,2} = 2\kappa\sigma^2 \left[\bm{\xi} (y-4x+3u_0)+ (By-4Ax+3u^*) \right]\\
\label{pre-xyzt16} &&\Sigma_{2,2} =  \kappa\sigma^2 \left[\bm{\xi} (z-4y+7x-4u_0)+(Cz-4By+7Ax-4u^*)\right]\\
\nonumber &&\Sigma_{3,2} =   \kappa\sigma^2 \left[\bm{\xi} (t-4z+6y-4x+u_0)+(Dt-4Cz+6By-4Ax+u^*)\right]\\
\nonumber &&
\\
\nonumber &&\Sigma_{1,3} =2\kappa\sigma^2\left[\bm{\zeta} (y-4x+3u_0) -2\alpha(by-2ax)\right]\\
\label{pre-xyzt17} &&\Sigma_{2,3} =  \kappa\sigma^2 \left[\bm{\zeta} (z-4y+7x-4u_0)-\alpha(3cz-8by+7ax)\right]\\
\nonumber &&\Sigma_{3,3} = \kappa\sigma^2\left[\bm{\zeta} (t-4z+6y-4x+u_0)-4\alpha(dt-3cz+3by-ax)\right]
\end{eqnarray}
where $\bm{\xi}=-\log \left(2\lambda^2\right)+\psi(H+1/2)$, $\bm{\zeta}=-2H/\lambda$, $u_0=g_H(0)$,
$u^*=\partial_H g_H(0)=(\log(2)+\psi(H))u_0$
 and
\begin{equation}
\label{xyzt}
\begin{cases}
\displaystyle
x= g_H(\lambda\alpha),~~ a=r_H(\lambda\alpha),& A=R_H(\lambda\alpha)
\\
\displaystyle
y= g_H(2\lambda\alpha),~~ b=r_H(2\lambda\alpha),& B=R_H(2\lambda\alpha)
\\
\displaystyle
z= g_H(3\lambda\alpha),~~ c=r_H(3\lambda\alpha),& C=R_H(3\lambda\alpha)
\\
\displaystyle
t= g_H(4\lambda\alpha),~~ d=r_H(4\lambda\alpha),& D=R_H(4\lambda\alpha),
\end{cases}
\end{equation}
with $r_H(\cdot)$ and $R_H(\cdot)$ being functions as defined in Lemma \ref{LEM0-App} and $\psi(\cdot)$ stands for the digamma function.
The next step is to introduce the concept of filters.
For a given $l,L\in\Z_+$ a filter of length $L+1$ and order $l$ is a sequence of reals $\AA=\left(a_0,a_1,\ldots,a_L\right)$ satisfying 
\begin{equation*}
\begin{cases}
\displaystyle
\sum_{q=0}^{L} a_q q^r =0, & \mbox{ for all } r=0,1,\ldots, (l-1)
\\
\displaystyle
\sum_{q=0}^{L} a_q q^l \neq 0,
\end{cases}
\end{equation*}with the convention $0^0=1$. The filter with zero order ($l=0$) is defined by $a_0=1$ and $a_q=0$ for all $0<q\leq L$.
Let $\AA_j$, $j=1,2,3$ be three filters of length $4$ defined as 
\begin{equation*}
\displaystyle
\AA_j:= \left(\AA_j(0),\AA_j(1),\AA_j(2),\AA_j(3)\right).
\end{equation*}
For each filter $\AA_j$ we associate a vector $\B_j = \left(\B_j(0),\B_j(1),\B_j(2),\B_j(3)\right)$, where
\begin{equation*}
\begin{cases}
\displaystyle
\B_j(0) = \sum_{q=0}^{3} \AA_j^2(q)
\\
\displaystyle
\B_j(k) = 2\sum_{q=0}^{3-k} \AA_j(k+q)\AA_j(q)
\end{cases}
\end{equation*}
for all $k = 1, 2, 3$. For practical puproses we choose  $\AA_1$ to be of order one, and  the others  $\AA_2$, $\AA_3$  to be of order two. Their numerical values are 
\begin{equation}
\label{eq11}
\begin{pmatrix}
\AA_1
\\
\AA_2
\\
\AA_3
\end{pmatrix}  = 
\displaystyle
\begin{pmatrix}
 0 & ~1 & ~0 & -1
\\
 0 & ~1 & -2 & ~1
\\
 1 & -1 & -1 & ~1
\end{pmatrix}.
\end{equation} 
After some elementary computations, we obtain the following associated $3\times 4$-matrix
\begin{equation}\label{eq12}
\displaystyle
\B =\begin{pmatrix}
\B_1
\\
\B_2
\\
\B_3
\end{pmatrix}  =\begin{pmatrix}
2 & ~0  & -2 & 0
\\
6 & -8 &  ~2 & 0
\\
4 & -2  & -4 & 2
\end{pmatrix}
\end{equation}
Let us assume that observations of the model \eqref{Model} at equidistant points $t_k=\alpha k$, $k=0,1,2,\ldots$ are available.
Using the method of generalized moments as proposed in \cite{barboza2017parameter},
we define our estimator $\widehat{\theta} = \left(\widehat{\sigma^2},\widehat{H},\widehat{\lambda}\right)$ by minimizing the contrast function
\begin{equation*}
\displaystyle
\widehat{\QQ}_N (\theta)= \sum_{l=1}^{3}\left(\frac {1}{N-2}\sum_{j=3}^N g_l \left(t_j,\theta \right)\right)^2,~~\theta=\left(\sigma^2,H,\lambda \right),
\end{equation*}
where $g_l(t,\theta)=\varphi_l^2(t)-V_l(\theta)$ and 
\begin{align*}
\displaystyle 
\begin{cases}
\displaystyle
\varphi_l(t):=\sum_{q=0}^{3} \AA_l(q) Y_\alpha (t-\alpha q),& \mbox{ for all }t\geq 0,
\\
\displaystyle
V_l(\theta) := \sum_{q=0}^{3} \B_l(q) \gamma_\theta (\alpha q),
\end{cases}
\end{align*}
where $\gamma_\theta(\cdot)$, $\AA_j$ and $\B$  are given respectively in \eqref{eq7.9}-\eqref{eq8}, \eqref{eq11} and  \eqref{eq12}.
The estimator $\widehat{\theta}$ is formally defined by
\begin{align}\label{MME}
 \widehat{\theta}_N &=\underset{\theta\in\Theta}{\arg\min}\widehat{\QQ}_N(\theta),
\end{align}

where $\Theta\subset \R^3$ is the parameter space. Note that this estimator depends on the step of observations $\alpha>0$ to be fixed. The following results are obtained.

\begin{theorem}
\label{TH2}
Consider the estimator defined by \eqref{MME} and assume that the true parameter $\theta = (\sigma^2,H,\lambda)$ belongs to the interior of the compact parameter space 
\begin{equation*}
\Theta \subset \left\lbrace (\sigma^2,H,\lambda)\in \R^3~:~\sigma^2>0,~H>1/2\mbox{ and }0<\lambda < e^{(\psi(H+1/2)+\psi(H))/2 } \right\rbrace
\end{equation*}with rectangular form, then  there exists some $\alpha^{*}>0$ depending only on $\Theta$ for which $\widehat{\theta}_N$ is a consistent estimator of $\theta$  almost surely and in $L^2(\Omega)$-sense, for every fixed  $\alpha\in(\alpha^*,\infty)$.
\end{theorem}
\begin{remark}\,\ \\ 
{\rm{\textbf{(i)}}} The range $(\alpha^{*},\infty)$ can be improved by adjusting $\delta\in(0,1)$ in \eqref{shConst} and optimizing the numerical constants ($M_j$, $L_j$, $j=0,1,2, \cdots$) appearing in the proof of Lemmas \ref{local-inj} and  \ref{principal-minors} to obtain  sharp bounds. The case where $H\in(0,1/2)$ can be treated by similar procedures. The main task is to verify the identifiability condition, specifically, to check whether the mapping $\theta\mapsto \varrho_{\theta,3}(\alpha)$ is globally injective (or univalent) on $\Theta$. For that purpose, one can rely on the existing global univalence results related to  (weak almost) P-matrices \cite{gale1965jacobian,parthasarathy1983global} or   (weak almost) N-matrices \cite{olech1991almost}, or PN-matrices \cite{inada1971production}. In general, due to the complexity of the functions involved in $\nabla \varrho_{\theta,3}(\alpha)$, it is not clear how to analytically verify such a condition  for a given $\alpha$. Nevertheless, this condition holds in the high-frequency estimation regime, that is, when the observation step    $\alpha$ is small ($\alpha\rightarrow 0$). This can be established by evaluating the principal minors given in Lemma \ref{principal-minors} by means of Taylor expansions.\\
 {\rm{\textbf{(ii)}}} In the case $H=1/2$ and $\lambda=0$, the process $Y_h$ reduces to the second order increments of the Brownian motion whose autocovariance is linear with respect to $\sigma^2$. Hence, the identifiability condition becomes trivial.
\end{remark}

\begin{proof}[Proof of Theorem \ref{TH2}] 
Here we adopt notations of \cite{barboza2017parameter}. Clearly, the matrix $\B$ defined by \eqref{eq12}
is of full rank $3$. There remains to verify the identifiability condition \cite[Assumption 2.1]{barboza2017parameter} for any $\alpha$ chosen in the range $(\alpha^{*},\infty)$, where $\alpha^{*}$ is the numerical constant due to Lemma \ref{principal-minors}. Let  $\UU\subset\R^3$ be an open rectangular region such that $\Theta \subset \UU $ and
\begin{equation}\label{constr}
\sigma^2>0,~H>1/2\mbox{ and }0<\lambda < e^{(\psi(H+1/2)+\psi(H))/2 } ,~\mbox{ for all }(\sigma^2,H,\lambda)\in\UU.
\end{equation} 
 Let $F:\UU\longrightarrow\R^3$ be the mapping defined by $F(\theta)=\Lambda \varrho_{\theta,3}(\alpha)$, where $\Lambda ={\rm{diag}}(-1,1,1)$ is a diagonal matrix. Clearly, $F$ is of class $C^1$ with the Jacobian matrix 
\begin{equation*}
\nabla F(\theta) =  \begin{pmatrix}
\displaystyle
-\Sigma_{1,1}~~ & -\Sigma_{1,2}~~ & -\Sigma_{1,3}
\\
\displaystyle
\Sigma_{2,1}~~ & \Sigma_{2,2}~~ & \Sigma_{2,3}
\\
\displaystyle
\Sigma_{3,1}~~ & \Sigma_{3,2}~~ & \Sigma_{3,3} 
\end{pmatrix}.
\end{equation*}
Taking into account that \eqref{constr} holds on $\UU$ and using Lemma \ref{principal-minors}, we assert that for any fixed $\alpha\in(\alpha^{*},\infty)$ the Jacobian matrix $\nabla F(\theta)$ has negative diagonals:
 $(\nabla F(\theta))_{1,1}=-\langle 2,3\rangle_\Sigma$, $(\nabla F(\theta))_{2,2}=\langle 1,3\rangle_\Sigma$ and $(\nabla F(\theta))_{3,3}=\langle 1,2\rangle_\Sigma$, while the remaining proper principal minors: $\langle 3\rangle_{\nabla F(\theta)}=-\langle 3\rangle_\Sigma$, $\langle 2\rangle_{\nabla F(\theta)}=-\langle 2\rangle_\Sigma $ and $\langle 1\rangle_{\nabla F(\theta)}=\langle 1\rangle_\Sigma  $ are negative. Furthermore, its determinant $\abs{\nabla F(\theta)}=-\langle\emptyset \rangle_\Sigma$ is positive. This means that  $\nabla F(\theta)$ is a weak almost N-matrix, and according to \cite[Theorem 4.1]{olech1991almost} with the associated remarks \cite[p. 123]{olech1991almost} $F$ is globally injective (or univalent) on $\UU$. Hence, the mapping $\theta\mapsto\varrho_{\theta,3}(\alpha)$ is globally injective on $\Theta$ as well.
 Moreover, the mapping $\theta\mapsto \B\varrho_{\theta}(\alpha)$ is injective (see Appendix B3).
  At this stage we can assert that all the conditions of \cite[Theorem 2.1 and Theorem 2.2]{barboza2017parameter} are met,
and the consistency of $\widehat{\theta}_N$ holds almost surely and in $L^2(\Omega)$-sense.
\end{proof}

\begin{theorem}
\label{TH3}
Under the assumptions of Theorem \ref{TH2} we have
\begin{equation*}
\displaystyle
\sqrt{N}\left(\widehat{\theta}_N-\theta\right) \overset{\LL}{\longrightarrow}
\NN \left(0,\CC(\theta)\Gamma (\theta) \CC (\theta)' \right)
\end{equation*}
as $N \to \infty$,   where
\begin{eqnarray*}
&&
\begin{cases}
\displaystyle
\left(\Gamma (\theta)\right)_{i,j} = 2\sum_{p\in\Z}\left(\sum_{k=0}^3 b^{i,j}_k \gamma_\theta \left(\alpha(p+k)\right)\right)^2,&~~\mbox{for all } i,j=1,2,3,
\\
\displaystyle
b^{i,j}_k := \sum_{q=0}^{3-k} \AA_i(q)\AA_j(q+k)+\sum_{q=0}^{3-k} \AA_j(q)\AA_i(q+k),&~~k=0,1,2,3
\end{cases}
\\ 
&&
\displaystyle
\CC(\theta)=- \left[\left(\nabla \varrho_{\theta}(\alpha) \right)'\B'\B \nabla\varrho_{\theta}(\alpha) \right]^{-1} \left(\nabla \varrho_{\theta}(\alpha) \right)'\B',
\end{eqnarray*}
$\nabla\varrho_{\theta}(\alpha) $ is the Jacobian of the vector $\varrho_{\theta}(\alpha):=\left(\gamma_\theta(0),\gamma_\theta(\alpha),\gamma_\theta(2\alpha),\gamma_\theta(3\alpha)\right)'$
with $\gamma_\theta(\cdot)$ being defined by \eqref{eq7.9}-\eqref{eq8} and $\AA_j$, $B$ are given in \eqref{eq11}-\eqref{eq12}.
\end{theorem}

\begin{proof}
We recall that the filters used in our framework are at most of order $2$ ($l_1=1$ and $l_2=l_3=2$).
According to \cite[Theorem 2.3]{barboza2017parameter} the estimator $\widehat{\theta}_N$
is asymptotically normal if the conditions of
Theorem \ref{TH2} are fulfilled and $\RR (\omega):=\left(1\wedge \vert \omega\vert^3\right)\overline{f}_\theta (\omega)\in L^2(-\pi,\pi)$, where 
\begin{equation*}
\displaystyle
\overline{f}_\theta (\omega) := \sum_{p\in\Z} f_\theta \left(\frac {\omega+2p\pi}{\alpha}\right),~~\omega\in [-\pi,\pi]
\end{equation*}
with $f_\theta(\cdot)$ being the spectral density function given in \eqref{sdf}
and $\theta$ denotes the true parameter.
Set $\Psi (\omega)=\left(1\wedge \vert \omega\vert^3\right)1_{(-\pi,\pi)}(\omega)$ and let $\norm{\Psi}$ be its $L^2$-norm.
Straightforward computations lead to 
\begin{eqnarray*}
&&
\displaystyle
\frac {\pi^2}{16\sigma^4\Gamma(H+1/2)^4}\int_{-\pi}^{\pi}\abs{\RR (\omega)}^2 d\omega
\\
&&
\displaystyle
\leq  \alpha^{4H+2}   \int \Psi^2(\omega)\left[\sum_{p\in\Z} \left((\lambda\alpha)^2+(\omega+2p\pi)^2\right)^{-H-1/2}\right]^2 d\omega
\\
&&
\displaystyle
\leq  3\alpha^{4H+2} \left\lbrace \int \frac { \Psi^2(\omega)}{\left((\lambda\alpha)^2+\omega^2\right)^{2H+1}} + \int \Psi^2(\omega)\left[\sum_{p\geq 1} \left((\lambda\alpha)^2+(\omega+2p\pi)^2\right)^{-H-1/2}\right]^2 d\omega \right.
\\
&&
\displaystyle
\left.~~~~~~~~~~~~~~ + \int \Psi^2(\omega)\left[\sum_{p\leq -1} \left((\lambda\alpha)^2+(\omega+2p\pi)^2\right)^{-H-1/2}\right]^2 d\omega \right\rbrace
\\
&&
\displaystyle
\leq  3\alpha^{4H+2}
\\
&&
\displaystyle
\cdot  \left\lbrace  \frac { \norm{\Psi}^2}{ (\lambda\alpha)^{4H+2}} +\frac {\left( \sum_{p\geq 1} \abs{2p-1}^{-2H-1}\right)^2 \norm{\Psi}^2}{\pi^{4H+2}}  +
\frac {\left( \sum_{p\leq 1} \abs{2p+1}^{-2H-1}\right)^2 \norm{\Psi}^2}{\pi^{4H+2}}\right\rbrace
\\
&&
\displaystyle
<\infty.
\end{eqnarray*}
This implies $\RR(\omega)\in L^2(-\pi,\pi)$.
In the second inequality above we split the series inside the integrand
as $\sum_{\lbrace p=0 \rbrace}+\sum_{p\geq 1}+\sum_{p\leq -1}$ and used the fact that $(u+v+w)^2\leq 3\left(u^2+v^2+w^2\right)$, for all $u,v,w\in\R$.
\end{proof}

\section{Alternative estimators of  $\sigma^2$, $H$ and $\lambda$}
\label{Sec.5}

As the drift $\mu$ can be estimated independently of the noise parameters,
we restrict ourselves to the model \eqref{Model} with $\mu=0$.
We denote the resulting process as $\bm{X}$, that is,
\begin{equation}
\label{Model'}
\displaystyle
\bm{X}_t=\sigma B_{H,\lambda}(t)
\end{equation}
for all $t \geq 0$.
Based on a discrete observations of  $\bm{X}$,
we aim to estimate $\sigma^2$, $H$ and $\lambda$.
Given a stochastic process $\xi=\lbrace \xi_t:t\geq 0 \rbrace$ and for $n\in\Z_+$ we write 
\begin{equation*}
\displaystyle
V_n^{(2)}(\xi):=\sum_{k=1}^{2^n} \left( \xi_{k/2^n}-\xi_{(k-1)/2^n}\right)^2
\end{equation*}
The following results are obtained.

\begin{theorem}\label{TH4}
Consider the model \eqref{Model'} with $H\in(0,1)$.
The following statements hold.
\begin{itemize}

\item[]\rm{\textbf{(i)}}
The statistic 
\begin{equation*}
\widetilde{H}_n = \frac{1}{2}\left(1-\frac{\log_2 V_n^{(2)}(\bm{X})}{n} \right) 
\end{equation*}
is a  strongly consistent estimator of $H$.
Furthermore, we have almost surely $\widetilde{H}_n-H=\OO(1/n)$, as $n\rightarrow\infty$. 

\item[]\rm{\textbf{(ii)}}
Under the condition $m/n=o(1)$, as $m,n\rightarrow\infty$,
the parameter $\sigma^2$ can be strongly estimated by
\begin{equation}\label{TILDsigma}
\widetilde{\sigma}_{m,n}^2=\frac {2\pi \sin (\pi \widetilde{H}_n)\Gamma(2\widetilde{H}_n +1)}{\Gamma(\widetilde{H}_n +1/2)^2} \cdot 2^{m(2\widetilde{H}_n-1)-1}V_m^{(2)}\left(\bm{X}\right)
\end{equation}  

\end{itemize}
\end{theorem}

\begin{remark}
If $H\in(0,1/2)$, then the same results of Theorem \ref{TH4} hold based on the observations of the model \eqref{Model}.
This is due to the fact that
\begin{equation*}
\displaystyle
2^{n(2H-1)}V_n^{(2)}(X) = 2^{n(2H-1)}V_n^{(2)} \left(\bm{X} \right)+o(1)
\end{equation*}
as $n \to \infty$ with probability one.
\end{remark}

\begin{proof}[Proof of Theorem \ref{TH4}]
It is well known that $\bm{X}$ is centered Gaussian with stationary increments, and admits the spectral density function 
\begin{eqnarray*}
\displaystyle
f_{\bm{X}} (\omega)
&=&
\displaystyle
\frac {\sigma^2\Gamma(H+1/2)^2}{2\pi \left(\lambda^2+\omega^2\right)^{H+1/2}}
\\
&=&
\displaystyle
\frac {\sigma^2\Gamma(H+1/2)^2}{2\pi\omega ^{2H+1}} \left(1-\OO \left(\omega^{-2}\right) \right)
\end{eqnarray*}
as $\omega \to \infty$.
So, we apply Lemma \ref{LEM4-App} with  $c= \sigma^2\Gamma(H+1/2)^2/(2\pi )>0$ and $\alpha=2H+1\in (1,3)$  to obtain the following almost sure convergences 
\begin{eqnarray}
\label{Q18A}
&&
\displaystyle
2-\frac {1}{n}\log_2 V_n^{(2)}(\bm{X})  \longrightarrow   \alpha,~~\mbox{as }n\rightarrow \infty,
\\
\label{Q19A}
&&
\displaystyle
-\cos \left(\frac {\alpha\pi}{2}\right) \Gamma(\alpha)2^{m(\alpha-2)-1}  V_m^{(2)} \left(\bm{X}\right)  \longrightarrow  c,~~\mbox{as }m\rightarrow \infty.
\end{eqnarray}
The first part of \rm{\textbf{(i)}} follows immediately from \eqref{Q18A} with the use of elementary computations.
For the second part, note that results of Lemma \ref{LEM4-App}
is just an application of \cite[Theorems 1 and 2]{gladyshev1961new} with $\gamma=3-\alpha$ and
$g(t)=-2c/\left( \cos\left(\frac {\alpha\pi}{2}\right)\Gamma(\alpha)\right)\neq 0$. So, revising the proof, one can readily derive a more precise statement showing the almost sure rate of convergence of $\widetilde{H}_n$, as given in  \rm{\textbf{(i)}}-Theorem \ref{TH4}.   The function $g(t)$ chosen here is well defined and nonnegative for $H\in(0,1)$. Clearly, this choice is not suitable for the boundary values of $H$. However, if $H=1$, one can still estimate $\sigma^2$ by a statistic with a reduced form. In fact, one has (see \cite[Theorems 4]{gladyshev1961new}) 
\begin{equation*}
\frac{2^{m+2}}{m\log(2)} V_m^{(2)}(\bm{X})  \overset{\P-a.s}{\longrightarrow} \sigma^2,~~\mbox{as }m\rightarrow \infty.
\end{equation*}

\rm{\textbf{(ii)}} 
First, equation \eqref{TILDsigma} can be rewritten as $\widetilde{\sigma}_{m,n}^2=\varphi_m (\widetilde{H}_n)V_m^{(2)}(\bm{X})$, where
\begin{equation}
\label{phi}
\displaystyle
\varphi_m(x) =\frac {2\pi \sin (\pi x)\Gamma(2x+1)}{\Gamma(x+1/2)^2} \cdot 2^{m(2x-1)-1},~~\mbox{ for all }x \in (0, 1).
\end{equation} 
 Performing some trigonometric computations with the use of \eqref{Q19A}, we can write 
\begin{eqnarray*}
\displaystyle
\widetilde{\sigma}_{m,n}^2
&=&
\displaystyle
\lbrace  \varphi_m(H) V_m^{(2)} (\bm{X})\rbrace\left(\frac {\varphi_m  (\widetilde{H}_n)}{\varphi_m(H)}\right)
\\
&\sim &
\displaystyle
\sigma^2 \left(\frac {\varphi_m (\widetilde{H}_n)}{\varphi_m(H)}\right),~~\P-a.s.
\end{eqnarray*}
To complete the proof, it suffices to prove that $\varphi_m (\widetilde{H}_n)/\varphi_m(H)\overset{\P-a.s}{\longrightarrow}1$
as $m,n\rightarrow\infty$.
By straightforward computations we have
\begin{equation}
\label{ratio-phi}
\displaystyle
\frac {\varphi_m(\widetilde{H}_n)}{\varphi_m(H)} = \frac {\sin (\pi\widetilde{H}_n)
\Gamma  (2\widetilde{H}_n+1 )}{\sin(\pi H)\Gamma(2H+1)}\cdot \left(\frac {\Gamma(H+1/2)}
{\Gamma (\widetilde{H}_n+1/2 )}\right)^2\cdot 2^{2m\left(\widetilde{H}_n-H\right)}.
\end{equation}
By means of continuous mapping theorem, both the first and the second
factors on the right hand side (RHS) of \eqref{ratio-phi} converge almost surely to $1$ as $n\rightarrow\infty$.
Now, using the last statement in \rm{\textbf{(i)}}-Theorem \ref{TH4} with
the hypothesis $m/n=o(1)$ yields $m (\widetilde{H}_n-H )=o(1)$ as $m,n\rightarrow\infty$.
Hence, the last factor on the RHS of \eqref{ratio-phi} converges as well to $1$, with probability one.
\end{proof}

Consider the estimators $\widetilde{H}$ and $\widetilde{\sigma}^2$ as defined in Theorem \ref{TH4} and the average 
\begin{equation*}
\displaystyle
\overline{\bm{X}}_n^{1,2}(h) = \frac {1}{n} \sum_{k=0}^{n-1}\left(\bm{X}_{h(k+1)}-\bm{X}_{hk}\right)^2,
\end{equation*}
then the following results are obtained.

\begin{theorem}
\label{TH5}
Consider the model \eqref{Model'} with $H\in(0,1)$.
The following statistic
\begin{equation*}
\displaystyle
\widetilde{\lambda}_{m,n,h} =\frac {1}{2}\exp \left(\frac {  \log \left(2\widetilde{\sigma}^2\Gamma (2\widetilde{H} )\right) -
\log\overline{\bm{X}}_n^{1,2}(h) }{2\widetilde{H}} \right)
\end{equation*}
converges in probability to $\lambda$ as $m,n,h\rightarrow\infty$.
\end{theorem}

\begin{proof}
First, let us prove that $\overline{\bm{X}}_n^{1,2}(h)\overset{L^2(\Omega)}{\longrightarrow} 2\sigma^2\Gamma(2H)/(2\lambda)^{2H} $ as $n,h\rightarrow\infty$.
For that purpose, we write
\begin{equation}
\label{Q22A}
\displaystyle
\overline{\bm{X}}_n^{1,2}(h) - \frac {2\sigma^2\Gamma(2H)}{(2\lambda)^{2H}} =
\underset{=:F_{n,h}}{\underbrace{\left(\overline{\bm{X}}_n^{1,2}(h)-\sigma^2\phi(h)\right)}} +
\sigma^2\left(\phi(h)-\frac {2\Gamma(2H)}{(2\lambda)^{2H}}\right),
\end{equation}
where $\phi(hx)=\abs{hx}^{2H}C^2_{hx}=:\phi_h(x)$.
In view of \eqref{VartfBm}, the last term on the RHS of \eqref{Q22A} converges towards $0$ as $h\rightarrow\infty$,
while the first term can be expressed as
\begin{eqnarray*}
\displaystyle
F_{n,h} = \frac {1}{n}\sum_{k=0}^{n-1}\left[\left(\widetilde{Z}_h(kh)\right)^2 - \E\left(\widetilde{Z}_{h}(kh)\right)^2\right],
\end{eqnarray*}
where $\widetilde{Z}_{h}(kh) =\bm{X}_{h(k+1)}-\bm{X}_{hk} $.
Clearly, $\E \left(F_{n,h}\right)=0$ and 
\begin{eqnarray}
\nonumber
&&
\displaystyle
\Var \left(F_{n,h}\right) =  n^{-2}\Var \left\lbrace \sum_{k=0}^{n-1} \left(\left(\widetilde{Z}_{h}(kh)\right)^2 - \E\left(\widetilde{Z}_{h}(kh)\right)^2\right) \right\rbrace
\\
\nonumber
&& =
\displaystyle
n^{-2}\Cov \left\lbrace \sum_{k=0}^{n-1}  \left(\left(\widetilde{Z}_{h}(kh)\right)^2-\E\left(\widetilde{Z}_{h}(kh)\right)^2\right),
\sum_{j=0}^{n-1}  \left(\left(\widetilde{Z}_{h}(jh)\right)^2-\E\left(\widetilde{Z}_{h}(jh)\right)^2\right)\right\rbrace
\\
\label{equality4} 
&& =
\displaystyle
n^{-2}\sum_{k=0}^{n-1}\sum_{j=0}^{n-1} \Cov     \left(\left(\widetilde{Z}_{h}(kh)\right)^2,
\left(\widetilde{Z}_{h}(jh)\right)^2 \right) = \frac {2}{n^2}\sum_{k=0}^{n-1}\sum_{j=0}^{n-1} \Cov^2     \left( \widetilde{Z}_{h}(kh),   \widetilde{Z}_{h}(jh)  \right)
\end{eqnarray}
\begin{eqnarray}
\nonumber
&&=
\displaystyle
\frac {2}{n^2}\sum_{k=0}^{n-1}\sum_{j=0}^{n-1} \Cov^2     \left( \widetilde{Z}_h \left(h(k-j)\right),
\widetilde{Z}_h(0)  \right),~~\mbox{(by stationarity of } \widetilde{Z}_{h}\mbox{)},
\\
\nonumber
&&=
\displaystyle
 \frac{\sigma^4}{2n^2}
 \sum_{k,j=0}^{n-1}\left( \phi_{h}(k-j+1)+\phi_{h}(k-j-1)-2\phi_{h}(k-j) \right)^2,~~\mbox{(by \eqref{eq4.9})},
\\
\nonumber
&& 
\displaystyle
=
  \frac{\sigma^4}{2n^2}
 \left(\sum_{k=j}+2\sum_{k=j-1}+2\sum_{k<j-1} \right)\left( \phi_{h}(k-j+1)+\phi_{h}(k-j-1)-2\phi_{h}(k-j) \right)^2
\\
\nonumber
&&
\displaystyle
\leq 
 \frac{\sigma^4}{2n^2}
  \left\lbrace 4n\phi_h(1)^2 +2n\left(\phi_h(2)-2\phi_h(1)\right)^2\right.
\\
\nonumber
&&
\displaystyle
\left.~~~~~~~+2\sum_{k<j-1} \left( \phi_{h}(k-j+1)+\phi_{h}(k-j-1)-2\phi_{h}(k-j) \right)^2 \right\rbrace
\\
\nonumber
&&
\displaystyle
\leq
 \frac{\sigma^4}{2n^2}
 \left\lbrace 4n\phi_h(1)^2 +2n \left(\phi_h(2)-2\phi_h(1)\right)^2\right.
\\
\label{lastSeries}
&& 
\displaystyle
\left. +2c_2^2\sum_{k<j-1} \left[2g_H \left(\lambda h\abs{k-j}\right)- g_H \left(\lambda h\abs{k-j+1}\right)-g_H \left(\lambda h\abs{k-j-1}\right) \right]^2 \right\rbrace, 
\end{eqnarray}
where $c_2=2\Gamma(H+1/2)/(\sqrt{\pi} 2^H \lambda^{2H})$ and $g_H(\cdot)$ being defined in the appendix.
In \eqref{equality4} we used the fact that  $\Cov \left(\xi^2,\zeta^2\right)=2\Cov^2(\xi,\zeta)$, 
for any Gaussian vector $(\xi,\zeta)$.
For the inequality \eqref{lastSeries}, we refer to \eqref{eq5}.
Now, by virtue of {\rm{\bf(i)}}-Lemma \ref{LEM(-1)-App}, the series appearing in \eqref{lastSeries} can be evaluated as 
\begin{eqnarray*}
&&
\displaystyle
\sum_{k<j-1} \left[2g_H \left(\lambda h\abs{k-j}\right) - g_H \left(\lambda h\abs{k-j+1}\right) - g_H\left(\lambda h\abs{k-j-1}\right) \right]^2
\\
&&
\displaystyle
\leq \AA_{\delta,H}^2  \sum_{j=2}^{n-1}\sum_{k=0}^{j-2}\left( 2e^{-\lambda \delta h\abs{k-j}}+
e^{-\lambda \delta h\abs{k-j+1}}+e^{-\lambda \delta h\abs{k-j-1}} \right)^2
\\
&&
\displaystyle
\leq 3 \AA_{\delta,H}^2  \sum_{j=2}^{n-1}\sum_{r=2}^{j}\left( 4e^{-2\lambda \delta hr}+e^{-2\lambda \delta h(r-1)}+e^{-2\lambda \delta h(r+1)} \right)
\leq \frac {18n \AA_{\delta,H}^2 }{1-e^{-2\lambda\delta h}}
\end{eqnarray*}
for some constant $\AA_{\delta,H}>0$ and $\delta\in(0,1)$.
Therefore, combining this with the fact that
\begin{eqnarray*}
\displaystyle
\lim_{h\rightarrow\infty}\phi_h(x)
&=&
\displaystyle
\lim_{h\rightarrow\infty}\phi(hx),~~\mbox{ for any  }x\neq 0
\\
&=&
\displaystyle
\lim_{h\rightarrow\infty}\abs{hx}^{2H} C^2_{hx}=2\Gamma(2H)/(2\lambda)^{2H},
\end{eqnarray*} 
the inequality \eqref{lastSeries} becomes 
\begin{eqnarray*}
\displaystyle
\Var \left(F_{n,h}\right) &&\leq 
\frac {2\sigma^4}{n} \left( 4\phi_h(1)^2 +2\left(\phi_h(2)-2\phi_h(1)\right)^2 +\frac {36c_2^2\AA_{\delta,H}^2 }{1-e^{-2\lambda \delta h}} \right)\\
&&\longrightarrow 0,
\end{eqnarray*}
as $h, n \to \infty$.
At this stage we can claim that $\overline{\bm{X}}_n^{1,2}(h)\overset{L^2(\Omega)}{\longrightarrow}2\sigma^2\Gamma(2H)/(2\lambda)^{2H} $
as $n,h\rightarrow\infty$.
Now, the continuous mapping theorem with the convergence of $\widetilde{H}$, $\widetilde{\sigma}^2$
and $\overline{\bm{X}}^{1,2}$ can be properly used to conclude the proof.
\end{proof}

\section{Conclusions}
\label{con}

In this paper, we have investigated the problem of
parameter estimation for a tempered fractional Brownian motion (tfBm)
with drift of the form $X_t = \mu t + \sigma B_{H,\lambda}(t)$, where $B_{H,\lambda}(t)$ denotes the tfBm.
By combining the method of moments with the ergodic theorem, we derived consistent and asymptotically normal estimators for the drift and noise parameters.

First, we showed that the drift estimator $\widehat{\mu} = X_{N}^{1, 1}(h)/h$ is unbiased,
strongly consistent, and asymptotically normal.
Next, by applying the generalized method of moments,
we constructed joint estimators for the parameters $\left(\sigma^2, H, \lambda\right)$,
proving their strong consistency and asymptotic normality under appropriate regularity conditions.
To address identifiability limitations, we adopted a two-step estimation strategy, separating the drift and noise estimation.

Additionally, using Gladyshev's limit theorems for Gaussian increment processes,
we proposed alternative, computationally efficient estimators for the parameters $\left( H, \sigma^2, \lambda \right)$
and established their strong or weak consistency.
These results extend existing methods for fractional and tempered fractional
models by providing a moment-based estimation framework that avoids heavy reliance on spectral or wavelet techniques.

Future research directions include the derivation of explicit asymptotic
distributions for the alternative estimators, the development of efficient numerical algorithms
for their implementation, and the extension of the proposed methodology to higher-order or multivariate tempered fractional models.

\section*{Appendix A. Technical lemmas}
\label{A}

In this section, we present some auxiliary lemmas needed to prove our main results.
Lemma \ref{LEM1-App} is well known  \cite[p. 413]{shiryaev1996probability},  while Lemma \ref{LEM(-1)-App} plays a vital role in our proofs, especially {\rm{\bf(i)}} which improves the upper bound given in \cite[p. 383]{gaunt2014inequalities}.
Consider the function $g_H$ defined on $(0,\infty)$
by $g_H(x)=x^H K_H(x)$, where $K_H(\cdot)$ is the modified Bessel function.
The norm in the space of bounded functions $f:(0,\infty)\longrightarrow\R$ is denoted as $\norm{f}_{0,\infty} =\sup_{x>0} \abs{f(x)}$.
The following results are obtained.
 
\begin{lemma}
\label{LEM(-1)-App}
For any $H>0$ the following statements hold.
\begin{itemize}
\item[]{\rm{\bf(i)}} For any $\delta\in(0,1)$, there exists a constant $\AA_{\delta,H}>0$ such that 
\begin{equation*}
\displaystyle
\abs{g_H(x)}\leq \AA_{\delta,H} e^{-\delta x}
\end{equation*}
for all $x > 0$. Moreover, $\AA_{\delta,H} \leq 2^{H-1}\Gamma(H)(1-\delta)^{-H}$.
\item[]{\rm{\bf(ii)}} We have $0< \frac{g_H(y)}{g_H(x)}\leq\left( \frac{y}{x}\right)^H e^{-(y-x)}$, for all $y\geq x>0$.
\end{itemize}
\end{lemma}
\begin{remark}\label{REM3}
The statement  {\rm{\bf(ii)}} is known since 1973  (e.g, \cite{baricz2010bounds} and references therein). We provide a different and simple proof of it. By using the identity $K_{-H}(x)=K_H(x)$ we deduce that {\rm{\bf(ii)}}-Lemma \ref{LEM(-1)-App} holds for $H<0$ as well.\\  

 If we combine {\rm{\bf(i)}} with the lower bound of $g_H$ given in \cite[p. 383]{gaunt2014inequalities}, we obtain 
\begin{equation}\label{eq1REM3}
2^{H-1}\Gamma(H)e^{-x} <g_H(x)\leq \AA_{\delta,H}e^{-\delta x},~~\mbox{for all } x>0.
\end{equation}If these inequalities are properly used with {\rm{\bf(ii)}}-Lemma \ref{LEM(-1)-App}, then for every $q>1$ one has
\begin{equation}\label{eq2REM3}
\frac{2^{H-1}\Gamma(H)}{\AA_{\delta,H}}e^{-(q-\delta) x} <\frac{g_H(q x)}{g_H( x)}< q^H e^{-(q-1) x},~~\mbox{for all } x>0.
\end{equation}
\end{remark}

\begin{proof}[Proof of Lemma \ref{LEM(-1)-App}]
{\rm{\bf(i)}} Let $p>1$, $a(p):=(p-1)/p$ and observe that  \cite[Eq 10.32.9]{olver2010nist}
\begin{equation*}
\displaystyle
K_H (x) =\int_0^\infty  e^{-x \cosh t} \cosh (Ht)dt
\end{equation*}
for all $H, x > 0$.
Therefore, for all $x>0$
\begin{eqnarray*}
\displaystyle
\abs{g_H(x)}
&=&
\displaystyle
\abs{x}^H \int_0^\infty  e^{-x a(p)\cosh t}\left( e^{-\frac {x}{p} \cosh t} \cosh (Ht) \right)dt
\\
&\leq &
\displaystyle
\abs{x}^H e^{-x a(p)} K_H\left(\frac {x}{p} \right) \leq  p^H \norm{g_H}_{0,\infty} e^{-x a(p)}\\
&\leq & 2^{H-1}p^H\Gamma(H)   e^{-x a(p)}.
\end{eqnarray*}In the last inequality above we used the fact that $\abs{g_H(x)}<2^{H-1}\Gamma(H)$, $\forall H,x>0$ (e.g., \cite[p. 383]{gaunt2014inequalities}). Due to the arbitrariness of $p$, we can choose its value so that $a(p)=\delta$.
Hence, the statement is achieved.     Once again, start from \cite[Eq 10.32.9]{olver2010nist} and write
\begin{eqnarray*}
\displaystyle
&& 0< g_H(y)= 
 \abs{y}^H \int_0^\infty  e^{-(y-x)\cosh t}\left( e^{-x \cosh t} \cosh (Ht) \right)dt\leq \abs{\frac{y}{x}}^He^{-(y-x)}g_H(x).
\end{eqnarray*}This proves {\rm{\bf(ii)}}, and the proof is complete. 
 
\end{proof}

\begin{lemma}
\label{LEM0-App}
For any $H>0$ the following statements hold true.
\begin{itemize}

\item[]{\rm{\textbf{(i)}}}
The partial derivatives of $g_H(z)$ are given as $\partial_H\lbrace g_H(z)\rbrace=R_H(z) g_H(z)$
and $\partial_z \lbrace g_H(z)\rbrace=-r_H(z) g_H(z)$, where $r_H(z):=\frac {K_{H-1}(z)}{K_H(z)}$ and 
\begin{align*}
\displaystyle
R_H(z):=\log(z) +2H\left(\int_z^\infty   \frac {I_H(s) K_H(s)}{s}ds -\frac {I_H(z)}{K_H(z)}\int_z^\infty  \frac {K_H^2(s)}{s}ds\right).
\end{align*}

\item[]{\rm{\textbf{(ii)}}}
$r_H(z)=1+o(1)$, as $z\rightarrow\infty$.

\item[]{\rm{\textbf{(iii)}}}
$R_H(z)=\log (z)+\frac {H}{2z}+o\left(\frac {1}{z}\right)$ as $z\rightarrow\infty$.

\item[]{\rm{\textbf{(iv)}}} $\abs{R_H(qz)-R_H(z)}\leq \log(q) +3M(2H+1) $, for all $q,z\in (1,\infty)$.

\item[]{\rm{\textbf{(v)}}} $R_H(z)\geq \log(z) - M$, for all $z\in (1,\infty)$.

\item[]{\rm{\textbf{(vi)}}} 
 $\vert R_H(z)\vert \leq \log(z) + M(6H+1)$, for all $z\in (1,\infty)$.
\end{itemize}
\end{lemma}

\begin{proof}
{\rm{\textbf{(i)}}}
It is well known that the modified Bessel function of the second kind $K_H(z)$
is differentiable with respect to $H$ and $z$.
Hence, $g_H(z)$ is differentiable as well.

Moreover, by virtue of \cite[(A.13)]{gaunt2014inequalities} one has $\partial_z g_H(z)=-z^H K_{H-1}(z)=-r_H(z) g_H(z)$.
On the other hand, taking derivatives with respect to $H$ with the use of \cite[(75)]{gonzalez2018closed} yields
\begin{eqnarray*}
&&
\displaystyle
\partial_H \lbrace g_H(z)\rbrace = \partial_H \lbrace z^H K_H(z) \rbrace
\\
&&=
\displaystyle
\log(z) z^H K_H(z)+2Hz^H\left( K_H(z)\int_z^\infty  \frac {I_H(s) K_H(s)}{s}ds -I_H(z)\int_z^\infty  \frac {K_H^2(s)}{s}ds \right)
\\
&&=
\displaystyle
g_H(z) \left[\log(z)+2H\left(\int_z^\infty  \frac {I_H(s) K_H(s)}{s}ds -\frac {I_H(z)}{ K_H(z)}\int_z^\infty  \frac {K_H^2(s)}{s}ds \right)\right]
\\
&&=:
\displaystyle
g_H(z)\left[\log(z)+2H \left( U_H(z)-\frac {I_H(z)}{ K_H(z)} V_H(z) \right)\right]
\\
&&=
\displaystyle
R_H(z) g_H(z)
\end{eqnarray*}
say.
The statement {\rm{\textbf{(ii)}}} is an immediate result of \cite[(A.5)]{gaunt2014inequalities}.
By using the  asymptotic expansions \cite[(67)-(68)]{gonzalez2018closed} of $I_H$ and $K_H$, as $z\rightarrow\infty$ we deduce that
\begin{eqnarray}
\label{IK}
&&
\displaystyle
K_H(z)^2\sim   \frac {\pi}{2}z^{-1}e^{-2z},~~  I_H(z)K_H(z)  \sim  \frac {1}{2z}
\\
\label{I/K}
&&
\displaystyle
\mbox{ and } \frac {I_H(z)}{K_H(z)}  \sim  \frac {e^{2z}}{\pi}
\end{eqnarray}
as $z \to \infty$.
Applying the L'Hopital's rule with the use of \eqref{IK} gives
$U_H(z)\sim \frac {1}{2z}$ and $V_H(z)\sim \frac {\pi}{4}\frac {e^{-2z}}{z}$ as $z\rightarrow\infty$.
Combine these results with \eqref{I/K} and perform some elementary computations to get the desired result {\rm{\textbf{(iii)}}}. Let  $q,z\in (1,\infty)$ then write the difference and rearrange the integrals as
\begin{align*}
 & R_H(qz)-R_H(z)=\log(z)+2H\left(\int_{qz}^\infty   \frac {I_H(s) K_H(s)}{s}ds -\frac {I_H(qz)}{K_H(qz)}\int_{qz}^\infty  \frac {K_H^2(s)}{s}ds\right)\\
 & -2H\left(\int_z^\infty   \frac {I_H(s) K_H(s)}{s}ds -\frac {I_H(z)}{K_H(z)}\int_z^\infty  \frac {K_H^2(s)}{s}ds\right)\\
 & =\log(q)-2H \int_z^{qz}   \frac {I_H(s) K_H(s)}{s}ds-2H \int_{qz}^\infty  \frac {K_H^2(s)}{s}ds\left(\frac {I_H(qz)}{K_H(qz)}-\frac {I_H(z)}{K_H(z)}\right)+2H\frac {I_H(z)}{K_H(z)}\int_z^{qz}\frac { K_H^2(s)}{s}ds \\
 &=: \log(q)-\JJ_1-\JJ_2+\JJ_3.
\end{align*}Applying \cite[(1.2)-Theorem 1]{baricz2016bounds} we claim that $I_H(s)K_H(s)\leq M s^{-1/3}$, for all $ s>0$. Here $M$ is the numerical constant defined as
\begin{equation}\label{M-const}
M~:=~\frac {2\pi^{3/2}} { \sqrt{3}\Gamma\left(\frac{2}{3}\right) \Gamma\left(\frac{5}{6}\right) \sqrt[3]{2} } \cdot \sup_{s\in\R_+}\sqrt[3]{s} J_0(s) 
\end{equation} and  $J_\nu$ stands for the Bessel function of the first kind. It follows that
\begin{eqnarray*}
   \abs{\JJ_1}  &=& 2H \int_z^{qz}   \frac {I_H(s) K_H(s)}{s}ds \leq 2MH \int_z^{qz}  s^{-4/3}  ds\\
&\leq & 6MH (1-q^{-1/3})z^{-1/3}<6MH .
\end{eqnarray*}and 
\begin{eqnarray}
 \nonumber  \abs{\JJ_2}  &=& 2H \int_{qz}^\infty  \frac {K_H^2(s)}{s}ds\abs{\frac {I_H(qz)}{K_H(qz)}-\frac {I_H(z)}{K_H(z)}}\\
  \nonumber  &\leq & 2H \left(\frac {I_H(qz)}{K_H(qz)}+\frac {I_H(z)}{K_H(z)}\right)\int_{qz}^\infty  \frac {g_H^2(s)}{s^{2H+1}}ds ,~~~~~~~(~g_H(s)=s^HK_H(s)~),\\
\label{eq1JJ2}    &\leq & 2H g_H^2(qz)\left(\frac {I_H(qz)}{K_H(qz)}+\frac {I_H(z)}{K_H(z)}\right)\int_{qz}^\infty   s^{-2H-1}ds \\
 \nonumber     &\leq & 2H g_H^2(qz)\left(\frac {I_H(qz)}{K_H(qz)}+\frac {I_H(z)}{K_H(z)}\right)\cdot \frac{(qz)^{-2H}}{2H}=K_H^2(qz)\left(\frac {I_H(qz)}{K_H(qz)}+\frac {I_H(z)}{K_H(z)}\right)\\
\nonumber      &\leq &  I_H(qz) K_H(qz)  +\left(\frac {K_H(qz)}{K_H(z)}\right)^2 [I_H(z)K_H(z)]\\
\nonumber      &\leq &  I_H(qz) K_H(qz)  +\left(\frac {g_H(qz)}{q^H g_H(z)}\right)^2 [I_H(z)K_H(z)]\\
\label{eq2JJ2}      &\leq & M (qz)^{-1/3}+q^{-2H} \cdot q^{2H} e^{-2(q-1)z}\cdot Mz^{-1/3}<2M.
\end{eqnarray}Equations \eqref{eq1JJ2} and \eqref{eq2JJ2} are justified by the fact that $g_H(\cdot)$ is decreasing on $(0,\infty)$ and \eqref{eq2REM3}, respectively. We proceed in a similar way to get
\begin{eqnarray}
\nonumber   \abs{\JJ_3}  &=& 2H\frac {I_H(z)}{K_H(z)}\int_z^{qz}\frac { K_H^2(s)}{s}ds  
 \nonumber     \leq   2H g_H^2(z)\frac {I_H(z)}{K_H(z)} \int_z^{qz} s^{-2H-1}ds \\
  \nonumber   &\leq  & 2H z^{2H}K_H^2(z)\frac {I_H(z)}{K_H(z)} \cdot \frac{(z^{-2H}-(qz)^{-2H})}{2H} =(1-q^{-2H})I_H(z)}{K_H(z)\\
 \label{eqJJ3}   &\leq & M (1-q^{-2H})z^{-1/3}<M.
\end{eqnarray} The statement {\rm{\textbf{(iv)}}} follows immediately by combining the previous results.
In view of \eqref{eqJJ3} we have 
\begin{align*}
& 2H\frac {I_H(z)}{K_H(z)}\int_z^{\infty}\frac { K_H^2(s)}{s}ds =\lim_{q\rightarrow \infty} 2H\frac {I_H(z)}{K_H(z)}\int_z^{qz}\frac { K_H^2(s)}{s}ds \leq M.
\end{align*}From this with $I_H(s),~K_H(s)>0$, $\forall s>0$ we derive the following 
\begin{eqnarray*}
R_H(z) &=& \log(z) +2H\left(\int_z^\infty   \frac {I_H(s) K_H(s)}{s}ds -\frac {I_H(z)}{K_H(z)}\int_z^\infty  \frac {K_H^2(s)}{s}ds\right)\\
&\geq & \log(z) -2H \frac {I_H(z)}{K_H(z)}\int_z^\infty  \frac {K_H^2(s)}{s}ds  \geq \log(z)-M,
\end{eqnarray*} which proves the last statement {\rm{\textbf{(v)}}}.
   As $I_H(z)$ and  $K_H(z)$ are positive functions, for $z>1$ we have 
\begin{eqnarray*}
&&  \vert R_H(z) \vert \leq \log(z) + \lim_{q\rightarrow\infty} 2H \int_z^{qz}   \frac {I_H(s) K_H(s)}{s}ds +\lim_{q\rightarrow\infty} 2H\frac {I_H(z)}{K_H(z)}\int_z^{qz}  \frac {K_H^2(s)}{s}ds.
\end{eqnarray*}Through the previous evaluation of $\vert\JJ_1\vert$ and $\vert\JJ_3\vert$ we conclude that {\rm{\textbf{(vi)}}} holds true. Therefore, the proof is complete.
\end{proof}

\begin{lemma}\rm{\textbf{(Birkhoff's ergodic theorem)}}
\label{LEM1-App}
Let $\lbrace  X_k\rbrace_{k=0}^\infty$ be a stationary ergodic sequence and let $g:\R\longrightarrow \R$
be a measurable function such that $\E \abs{g\left(X_0\right)}<\infty$, then 
\begin{equation*}
\displaystyle
\frac {1}{N} \sum_{k=0}^{N-1} g\left(X_k\right) \overset{\P-a.s}{\longrightarrow} \E g\left(X_0\right)
\end{equation*}
as $N \to \infty$.
\end{lemma}

Let us state the Gladyashev's theorem which is essential in establishing the asymptotics given Section \ref{Sec.5}. Consider a process with stationary Gaussian increments with  spectral decomposition of the form
\begin{equation}
Z(t)=Z(0)+\int_{-\infty}^\infty (e^{i\omega t}-1)d\widehat{B}(\omega),
\end{equation}where $\widehat{B}(\omega)$ is  a stochastic process with independent increments. If  $\E[\vert d\widehat{B}(\omega)\vert^2]=f_Z(\omega)d\omega$ then $f_Z(\omega)$ is said to be the spectral density of $Z$. Moreover, the following results are obtained (see \cite[p. 57]{gladyshev1961new})
\begin{lemma}\label{LEM4-App}
Let $Z(t)$ be a real stochastic process with stationary Gaussian increments and mean $0$, having the spectral density $f_Z(\omega)$  such that
\begin{equation}
\frac{f_Z(\omega)}{\vert\omega\vert^{-\alpha}}~=~c+\OO(\omega^{-2}),~~\mbox{as }\omega\rightarrow\infty,
\end{equation}
for some $c>0$ and $1<\alpha<3$. Then we have almost surely
\begin{equation}\label{Glad1}
 \lim_{n\rightarrow\infty} \left\lbrace2-\frac{1}{n}\log_2 \sum_{k=1}^{2^n}\left[Z\left(\frac{k}{2^n}\right)-Z\left(\frac{k-1}{2^n}\right)  \right]^2 \right\rbrace = \alpha\\
\end{equation} and 
\begin{equation}
\label{Glad2} \lim_{n\rightarrow\infty} \left\lbrace -\cos\left(\frac{\alpha\pi}{2}\right)\Gamma(\alpha)2^{n(\alpha-2)-1} \sum_{k=1}^{2^n}\left[Z\left(\frac{k}{2^n}\right)-Z\left(\frac{k-1}{2^n}\right) \right]^2 \right\rbrace = c
\end{equation}
\end{lemma}

\section*{Appendix B.}

\subsection*{B1. Local injectivity of the mapping $\theta \mapsto \varrho_{\theta,3}(\alpha)$}
\begin{lemma}\label{local-inj}
Let $\Theta \subset (0,\infty)\times (1/2,\infty) \times (0,\infty)$ be the compact parameter space, then there is some $\alpha_0>0$ depending only on $\Theta$ such that $\abs{\nabla \varrho_{\theta,3}(\alpha)}<0$, for all $\alpha\in(\alpha_0,\infty)$. In particular, the mapping $\theta \mapsto \varrho_{\theta,3}(\alpha)$ is locally injective.
\end{lemma}
\begin{proof}\label{proof-local-inj}
Consider the notations \eqref{eq9}, \eqref{xyzt} and set $\widetilde{\Sigma}_{1,1}=(y-4x+3u_0)$, $\widetilde{\Sigma}_{2,1}=(z-4y+7x-4u_0)$
and $\widetilde{\Sigma}_{3,1}=(t-4z+6y-4x+u_0)$.
Evaluating the determinant of $\nabla \varrho_{\theta,3}(\alpha)$ gives
\begin{align}
\nonumber  &
\displaystyle
\abs{\nabla \varrho_{\theta,3}(\alpha)}\\
\nonumber  &= 2\kappa^3\sigma^4
\cdot  \begin{vmatrix}
\displaystyle
\widetilde{\Sigma}_{1,1} &~~ \bm{\xi} \widetilde{\Sigma}_{1,1}+ (By-4Ax+3u^*) &~~ \bm{\zeta} \widetilde{\Sigma}_{1,1} -2\alpha(by-2ax)
\\
\displaystyle
\widetilde{\Sigma}_{2,1} &~~ \bm{\xi} \widetilde{\Sigma}_{2,1}+(Cz-4By+7Ax-4u^*) &~~ \bm{\zeta} \widetilde{\Sigma}_{2,1}-\alpha(3cz-8by+7ax)
\\
\displaystyle
\widetilde{\Sigma}_{3,1} &~~ \bm{\xi} \widetilde{\Sigma}_{3,1}+(Dt-4Cz+6By-4Ax+u^*) &~~ \bm{\zeta} \widetilde{\Sigma}_{3,1}-4\alpha(dt-3cz+3by-ax)
\end{vmatrix}
\\
\nonumber &
= 2\kappa^3\sigma^4 \begin{vmatrix}
\displaystyle
\widetilde{\Sigma}_{1,1} &~~   (By-4Ax+3u^*) &~~   -2\alpha(by-2ax)
\\
\displaystyle
\widetilde{\Sigma}_{2,1} &~~  (Cz-4By+7Ax-4u^*) &~~  -\alpha(3cz-8by+7ax)
\\
\displaystyle
\widetilde{\Sigma}_{3,1} &~~(Dt-4Cz+6By-4Ax+u^*) &~~  -4\alpha(dt-3cz+3by-ax)
\end{vmatrix}
\\
\nonumber &
= -2\alpha\kappa^3\sigma^4 \begin{vmatrix}
\displaystyle
\widetilde{\Sigma}_{1,1} &~~   (By-4Ax+3u^*) &~~    2 (by-2ax)
\\
\displaystyle
\widetilde{\Sigma}_{2,1} &~~  (Cz-4By+7Ax-4u^*) &~~  (3cz-8by+7ax)
\\
\displaystyle
\widetilde{\Sigma}_{3,1} &~~(Dt-4Cz+6By-4Ax+u^*) &~~  4(dt-3cz+3by-ax)
\end{vmatrix}
\end{align}
\begin{align}
\nonumber &=
-2\alpha\kappa^3\sigma^4
\\
\nonumber &
\cdot \begin{vmatrix}
\displaystyle
\widetilde{\Sigma}_{1,1} &~~   (B-A)y-3A^* &~~    (2b-a)y-3au_0
\\
\displaystyle
\widetilde{\Sigma}_{2,1} &~~  (C-A)z+4(A-B)y+4A^* &~~  (3c-a)z+4(a-2b)y+4au_0
\\
\displaystyle
\widetilde{\Sigma}_{3,1} &~~(D-A)t+4(A-C)z+6(B-A)y-A^* &~~  (4d-a)t+4(a-3c)z+6(2b-a)y-au_0
\end{vmatrix}
\\
& ~~~~~~~~~~~~~~~~~~~~~~~~~~~~~~~~~~~~~~~~~~~~~~~~~~~~~~~~~A^*=Au_0-u^*
\\
\nonumber &=
-2\alpha\kappa^3\sigma^4
\\
\nonumber &
\cdot \begin{vmatrix}
\displaystyle
(y-4x)+3u_0 &~~   (B-A)y-3A^* &~~    (2b-a)y-3au_0
\\
\displaystyle
(z-9x)+8u_0 &~~  (C-A)z-8A^*  &~~  (3c-a)z-8au_0
\\
\displaystyle
\widetilde{\Sigma}_{3,1} &(D-A)t+4(A-C)z+6(B-A)y-A^* &  (4d-a)t+4(a-3c)z+6(2b-a)y-au_0
\end{vmatrix}
\\
\label{InterStep} &=
-2\alpha\kappa^3\sigma^4  \begin{vmatrix}
\displaystyle
(y-4x)+3u_0 &~~   (B-A)y-3A^* &~~   (2b-a)y-3au_0
\\
\displaystyle
(z-9x)+8u_0 &~~  (C-A)z-8A^*  &~~  (3c-a)z-8au_0
\\
\displaystyle
(t-16x)+15u_0 &~~(D-A)t-15A^*  &~~  (4d-a)t-15au_0
\end{vmatrix} \\
\nonumber &= -2\alpha\kappa^3\sigma^4 \left\lbrace \begin{vmatrix}
\displaystyle
(y-4x)  &~~   (B-A)y  &~~   (2b-a)y 
\\
\displaystyle
(z-9x)  &~~  (C-A)z   &~~  (3c-a)z 
\\
\displaystyle
(t-16x)  &~~(D-A)t   &~~  (4d-a)t 
\end{vmatrix}+\begin{vmatrix}
\displaystyle
(y-4x)  &~~   (B-A)y  &~~    -3au_0
\\
\displaystyle
(z-9x)  &~~  (C-A)z   &~~   -8au_0
\\
\displaystyle
(t-16x)  &~~(D-A)t   &~~   -15au_0
\end{vmatrix}\right.\\
\nonumber & \left. +\begin{vmatrix}
\displaystyle
(y-4x)  &~~    -3A^* &~~   (2b-a)y 
\\
\displaystyle
(z-9x)  &~~   -8A^*  &~~  (3c-a)z 
\\
\displaystyle
(t-16x)  &~~ -15A^*  &~~  (4d-a)t 
\end{vmatrix} +\begin{vmatrix}
\displaystyle
 3u_0 &~~   (B-A)y  &~~   (2b-a)y 
\\
\displaystyle
 8u_0 &~~  (C-A)z   &~~  (3c-a)z 
\\
\displaystyle
 15u_0 &~~(D-A)t   &~~  (4d-a)t 
\end{vmatrix}\right\rbrace 
\end{align}
\begin{align}
\nonumber &= -2\alpha\kappa^3\sigma^4 \left\lbrace xy\begin{vmatrix}
\displaystyle
(y/x-4)  &   (B-A)  &    (2b-a)y 
\\
\displaystyle
(z/x-9)  &   (C-A)z/y   &   (3c-a)z 
\\
\displaystyle
(t/x-16)  & (D-A)t/y   &   (4d-a)t 
\end{vmatrix}+xy\begin{vmatrix}
\displaystyle
(y/x-4)  &    (B-A)   &     -3au_0
\\
\displaystyle
(z/x-9)  &   (C-A)z/y   &    -8au_0
\\
\displaystyle
(t/x-16)  & (D-A)t/y   &    -15au_0
\end{vmatrix}\right.\\
\nonumber & \left. +A^*xy\begin{vmatrix}
\displaystyle
(y/x-4)  &     -3 &    (2b-a)  
\\
\displaystyle
(z/x-9)  &    -8   &   (3c-a)z/y 
\\
\displaystyle
(t/x-16)  &  -15   &   (4d-a)t/y 
\end{vmatrix} +xy\begin{vmatrix}
\displaystyle
 3u_0 &    (B-A)y/x  &    (2b-a) 
\\
\displaystyle
 8u_0 &   (C-A)z/x   &   (3c-a)z/y 
\\
\displaystyle
 15u_0 & (D-A)t/x   &   (4d-a)t/y 
\end{vmatrix}\right\rbrace 
\\
\nonumber &= -2\alpha\kappa^3\sigma^4 \left\lbrace
o(xy)+\OO(xy) +87A^*xy(1+o(1)) +o(xy)
\right\rbrace,~~~~\mbox{as }\alpha\rightarrow \infty \\
\nonumber &\sim  -2\alpha\kappa^3\sigma^4A^*xy,~~~~\mbox{as }\alpha\rightarrow \infty.
\end{align}We justify the above asymptotics by means of \eqref{eq1REM3}, \eqref{eq2REM3} and Lemma \ref{LEM0-App}. For example $y/x=g_H(2\lambda\alpha)/g_H(\lambda\alpha)=\OO(e^{-\lambda\alpha})$, $a=r_H(\lambda\alpha)=1+o(1)$, $b=r_H(2\lambda\alpha)=1+o(1)$  and  $A^*=Au_0-u^*=u_0R_H(\lambda\alpha)-u^*\sim u_0\log(\lambda\alpha)$, etc, as $\alpha\rightarrow \infty$. Although $\abs{\nabla \varrho_{\theta,3}(\alpha)}$ converges rapidly to zero, as $\alpha\rightarrow\infty$, we can find some $\alpha$ large enough so that $\abs{\nabla \varrho_{\theta,3}(\alpha)}<0$. This proves the identifiability condition for $\alpha$ large enough. For a practical purpose, we need to specify the range of $\alpha$, for which  $\abs{\nabla \varrho_{\theta,3}(\alpha)}< 0$.\bigskip

Observe that each column in \eqref{InterStep} can be expressed as $\left(F(q,\xi)\right)_{q=2,3,4}$ with $\xi=\lambda\alpha$. So,   we simplify notations by denoting the determinant as $ {\bf\rm{det}} \left\lbrace  \cdot ;\cdot;\cdot \right\rbrace$ where each expression inside it  represents a column vector. The following additional notations are considered
\begin{eqnarray*}
&& g(\xi)=g_H(\xi),~~g'(\xi)=g_{H-1}(\xi),~~R(\xi)=R_H(\xi),\\
&& \varphi_q(\xi)=\frac{g(q\xi)}{g(\xi)},~~\varphi_q'(\xi)=\frac{g'(q\xi)}{g'(\xi)}\\
&& x'=g'(\xi),~~ y'=g'(2\xi),~~ z'=g'(3\xi),~~ t'=g'(4\xi).
\end{eqnarray*} 
 By definition (see Lemma \ref{LEM0-App}) of $r_H(\cdot)$ we have $a=r_H(\xi)=\xi(x'/x)$. Similarly, $b=2\xi(y'/y)$, $c=3\xi(z'/z)$ and $d=4\xi(t'/t)$. Therefore, equation \eqref{InterStep} becomes 
\begin{eqnarray}
\nonumber && \abs{\nabla \varrho_{\theta,3}(\alpha)} =-2\alpha\kappa^3\sigma^4  \begin{vmatrix}
\displaystyle
(y-4x)+3u_0 &~~   (B-A)y-3A^* &~~   \left(4\xi\frac{y'}{y}-\xi\frac{x'}{x}\right)y-3u_0\xi\frac{x'}{x}
\\
\displaystyle
(z-9x)+8u_0 &~~  (C-A)z-8A^*  &~~  \left(9\xi\frac{z'}{z}-\xi\frac{x'}{x}\right)z-8 u_0\xi\frac{x'}{x}
\\
\displaystyle
(t-16x)+15u_0 &~~(D-A)t-15A^*  &~~  \left(16\xi\frac{y'}{y}-\xi\frac{x'}{x}\right)t-15 u_0\xi\frac{x'}{x}
\end{vmatrix}\\
\nonumber &&=-2\alpha\kappa^3\sigma^4 \frac{\xi}{x} \begin{vmatrix}
\displaystyle
(y-4x)+3u_0 &~~   (B-A)y-3A^* &~~    (4xy'-x'y) -3u_0x'
\\
\displaystyle
(z-9x)+8u_0 &~~  (C-A)z-8A^*  &~~   (9xz'-x'z) -8 u_0x'
\\
\displaystyle
(t-16x)+15u_0 &~~(D-A)t-15A^*  &~~  (16xt'-x't) -15 u_0x'
\end{vmatrix}
\\
\nonumber &&=-2\alpha\kappa^3\sigma^4 \frac{\xi}{x} \begin{vmatrix}
\displaystyle
(y-4x)+3u_0 &~~   (B-A)y-3A^* &~~    4(xy'-xx') 
\\
\displaystyle
(z-9x)+8u_0 &~~  (C-A)z-8A^*  &~~  9(xz'-xx')  
\\
\displaystyle
(t-16x)+15u_0 &~~(D-A)t-15A^*  &~~   16(xt'-xx')  
\end{vmatrix}
\\
\nonumber &&=-2\alpha\kappa^3\sigma^4  \xi  \begin{vmatrix}
\displaystyle
(y-4x)+3u_0 &~~   (B-A)y-3A^* &~~    4(y'-x') 
\\
\displaystyle
(z-9x)+8u_0 &~~  (C-A)z-8A^*  &~~  9(z'-x')  
\\
\displaystyle
(t-16x)+15u_0 &~~(D-A)t-15A^*  &~~   16(t'-x')  
\end{vmatrix}\\
\nonumber &&= -2\alpha\kappa^3\sigma^4  \xi {\bf\rm{det}} \left\lbrace
[g(q\xi)-q^2g(\xi)]+(q^2-1)u_0;    [R(q\xi)-R(\xi)]g(q\xi)+(1-q^2)A^*;     q^2\left[ g'(q\xi)-g'(\xi)\right] 
\right\rbrace\\
\nonumber && 
=   -2\alpha\kappa^3\sigma^4  \xi g(2\xi)g'(\xi) A^* {\bf\rm{det}}\left\lbrace \CC_1+\CC_2;\CC_3+\CC_4;\CC_5\right\rbrace
\\
\label{DET} && 
= -2\alpha\kappa^3\sigma^4  \xi g(2\xi)g'(\xi) A^* \left[ {\bf\rm{det}} \left\lbrace\CC_1;\CC_3;\CC_5\right\rbrace+{\bf\rm{det}} \left\lbrace\CC_1;\CC_4;\CC_5\right\rbrace+{\bf\rm{det}} \left\lbrace\CC_2;\CC_3;\CC_5\right\rbrace \right]
\end{eqnarray}where 
\begin{eqnarray}
\label{eq1DETnots} && A^*=u_0R(\xi)-u^*,u_0=2^{H-1}\Gamma(H)~~
 \CC_1=g(\xi) [\varphi_q(\xi)-q^2],
~~\CC_2=(q^2-1)u_0,\\
\label{eq2DETnots} && \CC_3=\frac{g(q\xi)}{A^*}[R(q\xi)-R(\xi)],~~\CC_4=(1-q^2),~~
\CC_5=[ q^2\varphi_q'(\xi)-q^2]/g(2\xi)
\end{eqnarray}
Assume that   $H\in\KK:=[\underline{H},\overline{H}]$ and $\lambda\in [\underline{\lambda},\overline{\lambda}]$. This hypothesis is technical and used to define the constants $M^*$, $M_j$, $j=0,1,\cdots,5$ by maximizing functions in $H$ or $\lambda$. So, it can be dropped by considering the same functions in $\theta=(\sigma^2,H,\lambda)$ and maximizing them over the parameter space $\Theta$. Let $\delta\in(0,1)$ to be specified later. We shall start the evaluation of the determinants appearing in \eqref{DET}.
\paragraph*{\textbf{i-Evaluating the determinant $\bm{{\bf\rm{det}} \left\lbrace\CC_1;\CC_3;\CC_5\right\rbrace}$ .}}
 We apply \eqref{eq1REM3} and \eqref{eq2REM3} to get
\begin{equation}\label{eq1INQ's}
\begin{cases}
 u_0 e^{-q\xi}
  <g(q\xi)\leq u_0(1-\delta)^{-H}e^{-q\delta \xi},&\\
  (1-\delta)^{H}e^{-(q-\delta) \xi} <\varphi_q(\xi)< q^H e^{-(q-1) \xi},&\\
    ~&\\
 0< \varphi_q'(\xi)< q^{H-1} e^{-(q-1) \xi} & 
\end{cases} 
\end{equation}for all $\xi,H >0$ and $q=2,3,4$. In particular, for all $\xi>0$ 
\begin{equation}\label{Ineq-varphi2}
1/\varphi_2(\xi)\leq (1-\delta)^{-H} e^{(2-\delta)\xi} 
\end{equation} and 
\begin{equation}\label{Ineq-varphiq-2}
\begin{cases}
g(q\xi)/g(2\xi) =1,~~\mbox{ for }q=2\\
 g(q\xi)/g(2\xi) =\varphi_q(\xi)/\varphi_2(\xi)\leq q^H(1-\delta)^{-H}e^{-[q-3+\delta]\xi},~~\mbox{ for }q=3,4.
\end{cases}
\end{equation}
Keep in mind that the inequalities \eqref{eq1INQ's}-\eqref{Ineq-varphiq-2} will be systematically used.  Straightforward computations with the use of Lemma \ref{LEM0-App} gives  
\begin{align}\label{eq2INQ's}
\nonumber \abs{\frac{\CC_{3}}{g(2\xi)}} &=\frac{1}{\vert A^*\vert}\abs{[R(q\xi)-R(\xi)]} \cdot\frac{g(q\xi)}{g(2\xi)}\\
\nonumber &\leq  \frac{1}{\vert A^*\vert}\left[\log(4)+3M(2H+1)\right]4^H(1-\delta)^{-H} \\
&\leq \frac{M_{0}}{\vert A^*\vert},~~\mbox{for all } \xi>1,~\mbox{and }q=2,3,4.\\
\nonumber \mbox{where } &   M_{0}=\sup_{H\in\KK} \left\lbrace 4^H(1-\delta)^{-H}  [\log(4) +3M(2H+1)]\right\rbrace~\mbox{and }M\mbox{ being defined by \eqref{M-const}}.
\end{align} But we know that 
\begin{eqnarray*}
&& {\bf\rm{det}} \left\lbrace\CC_{1};\CC_{3};\CC_{5}\right\rbrace =
 {\bf\rm{det}} \left\lbrace g(\xi)[\varphi_q(\xi)-q^2];\CC_{3};[ q^2\varphi_q'(\xi)-q^2]/g(2\xi)\right\rbrace\\
&& = g(\xi)\left[ {\bf\rm{det}} \left\lbrace \varphi_q(\xi);\frac{\CC_{3}}{g(2\xi)};  q^2\varphi_q'(\xi) \right\rbrace +  {\bf\rm{det}} \left\lbrace  \varphi_q(\xi) ;\frac{\CC_{3}}{g(2\xi)}; -q^2 \right\rbrace 
+  {\bf\rm{det}} \left\lbrace -q^2;\frac{\CC_{3}}{g(2\xi)};  q^2\varphi_q'(\xi) \right\rbrace \right].
\end{eqnarray*}
Throughout the evaluation of this determinant at least one of  the terms: $\varphi_q(\xi)$ and $\varphi'_q(\xi)$ arises with $\CC_{3}/g(2\xi)$. So, we use \eqref{eq1INQ's} and \eqref{eq2INQ's} to arrive at

\begin{align}
\nonumber  \abs{{\bf\rm{det}} \left\lbrace\CC_{1};\CC_{3};\CC_{5}\right\rbrace} &\leq 16\cdot 3! M_{0}\cdot\frac{g(\xi)}{\vert A^*\vert}\left( 4^{2H-1} e^{-2\xi}+4^{H} e^{-\xi}+16\cdot 4^{H-1} e^{-\xi}  \right)
 \\
\nonumber &\leq 96u_0(1-\delta)^{-H} M_{0}\cdot\frac{e^{-(1+\delta)\xi}}{\vert A^*\vert}\left( 4^{2H-1} e^{-\xi}+5\cdot 4^{H} \right)
\\
\label{eq3INQ's} & \leq \frac{M_{1}}{\vert A^*\vert} e^{-(1+\delta)\xi},
\end{align}where $M_1:=96 M_0\sup_{H\in\KK}\left\lbrace u_0(1-\delta)^{-H}( 4^{2H-1}  +5\cdot 4^{H} ) \right\rbrace$.

\paragraph*{\textbf{ii-Evaluating the determinant $\bm{{\bf\rm{det}} \left\lbrace\CC_2;\CC_3;\CC_5\right\rbrace}$ .}}
In similar fashion we obtain
\begin{eqnarray}
\nonumber &&  {\bf\rm{det}} \left\lbrace\CC_2;\CC_{3};\CC_{5}\right\rbrace  = {\bf\rm{det}} \left\lbrace (q^2-1)u_0;\CC_{3};[ q^2\varphi_q'(\xi)-q^2]/g(2\xi)\right\rbrace
\\
\nonumber  &&= u_0\left(  {\bf\rm{det}} \left\lbrace (q^2-1);\frac{\CC_{3}}{g(2\xi)};  q^2\varphi_q'(\xi) \right\rbrace +   {\bf\rm{det}} \left\lbrace  (q^2-1) ;\frac{\CC_{3}}{g(2\xi)}; -q^2 \right\rbrace\right).
\end{eqnarray}
 Thus,
\begin{eqnarray}
 \nonumber && \abs{ {\bf\rm{det}} \left\lbrace\CC_2;\CC_{3};\CC_{5}\right\rbrace }  \leq 16\times 15\cdot 3! u_0\cdot\frac{M_{0}}{\vert A^*\vert}\left( 4^{H-1}e^{-\xi}+1\right) \\  
\label{eq4INQ's} && \leq M_2\frac{u_0}{\vert A^*\vert} 
\end{eqnarray}where $M_2:=1440  M_{0}\sup_{H\in\KK} \lbrace ( 4^{H-1}+1) \rbrace$.
Now, we recall that $u^*=u_0(\log(2)+\psi(H))$, $u_0=2^{H-1}\Gamma(H)$ and $A=R(\xi)$. From {\rm{\textbf{(v)}}}-Lemma \ref{LEM0-App} we deduce that
\begin{align*}
\nonumber A^* =u_0R(\xi)-u^*  &\geq u_0 \log(\xi)-u_0[M+\log(2)+\psi(H)],~~\mbox{(see \eqref{M-const})}\\
  &\geq u_0 [ \log(\xi)-M^*],~~\mbox{ for all }\xi > e^{M^*},
\end{align*}with $M^*=\sup_{H\in\KK}\left\lbrace M+\log(2)+\psi(H) +1/u_0\right\rbrace$.
 Consequently, from \eqref{eq4INQ's} we derive
\begin{equation}\label{eq8INQ's}
\abs{{\bf\rm{det}} \left\lbrace\CC_2;\CC_{3};\CC_{5}\right\rbrace} \leq \frac{M_2}{ \log(\xi)-M^*},~~\mbox{ for all }\xi > 1\vee e^{M^*}.
\end{equation}Note also that by definition of the constant $M^*$, we have
\begin{equation}\label{eq7INQ's}
A^*\geq 1,~~\mbox{ for all }\xi > 1\vee e^{M^*}.
\end{equation}

\paragraph*{\textbf{iii-Evaluating the determinant $\bm{{\bf\rm{det}} \left\lbrace\CC_1;\CC_4;\CC_5\right\rbrace}$ .}}
Once again, we proceed in similar fashion and write
\begin{eqnarray*}
&& {\bf\rm{det}} \left\lbrace\CC_{1};\CC_{4};\CC_{5}\right\rbrace =  {\bf\rm{det}} \left\lbrace g(\xi) [\varphi_q(\xi)-q^2];(1-q^2);[ q^2\varphi_q'(\xi)-q^2]/g(2\xi)\right\rbrace  
\\
&& =  {\bf\rm{det}} \left\lbrace  [\varphi_q(\xi)-q^2];(1-q^2);[ q^2\varphi_q'(\xi)-q^2]\right\rbrace  
\cdot [g(\xi)/g(2\xi)]
\\
&& =\left[ {\bf\rm{det}} \left\lbrace  \varphi_q(\xi) ;(1-q^2);  q^2\varphi_q'(\xi) \right\rbrace
+{\bf\rm{det}} \left\lbrace  \varphi_q(\xi) ;(1-q^2); -q^2 \right\rbrace
+{\bf\rm{det}} \left\lbrace  -q^2 ;(1-q^2);  q^2\varphi_q'(\xi) \right\rbrace \right]
/\varphi_2(\xi)
\\
&& = \left[ {\bf\rm{det}} \left\lbrace  \varphi_q(\xi) ;(1-q^2);  q^2\varphi_q'(\xi) \right\rbrace
+{\bf\rm{det}} \left\lbrace  \varphi_q(\xi) ;1; -q^2 \right\rbrace
+{\bf\rm{det}} \left\lbrace  -q^2 ;1;  q^2\varphi_q'(\xi) \right\rbrace \right]
/\varphi_2(\xi)
\\
&& = \left[{\bf\rm{det}} \left\lbrace  \varphi_q(\xi) ;(1-q^2);  q^2\varphi_q'(\xi) \right\rbrace
-{\bf\rm{det}} \left\lbrace  \varphi_q(\xi) ;1; q^2 \right\rbrace
+{\bf\rm{det}} \left\lbrace q^2\varphi_q'(\xi)  ;1;  q^2  \right\rbrace \right]
/\varphi_2(\xi)
\\
&& =\left[ {\bf\rm{det}} \left\lbrace  \varphi_q(\xi) ;(1-q^2);  q^2\varphi_q'(\xi) \right\rbrace
+{\bf\rm{det}} \left\lbrace \Phi_q  ;1;  q^2  \right\rbrace \right]
/\varphi_2(\xi),
\end{eqnarray*}where  $\Phi_q=q^2\varphi_q'(\xi) -\varphi_q(\xi)$. The last determinant simplifies to 
\begin{eqnarray*}
\nonumber  {\bf\rm{det}}\left\lbrace \Phi_q  ;1;  q^2  \right\rbrace 
 &=& \begin{vmatrix}
\Phi_2  & 1 &   4\\
\Phi_3  & 1 &   9\\
\Phi_4  & 1 &   16
\end{vmatrix}  =\Phi_2\begin{vmatrix}
 1 &   9\\
 1 &   16
\end{vmatrix}-\Phi_3\begin{vmatrix}
   1 &   4\\
  1 &   16
\end{vmatrix}+\Phi_4\begin{vmatrix}
  1 &   4\\
  1 &   9
\end{vmatrix}\\
&=& 7\Phi_2-12\Phi_3+5\Phi_4.
\end{eqnarray*} Thereby, 
\begin{eqnarray}
\nonumber && {\bf\rm{det}} \left\lbrace\CC_{1};\CC_{4};\CC_{5}\right\rbrace = 7 \Phi_2/\varphi_2(\xi)+ {\bf\rm{det}} \label{eq7INQ's-bis} \left\lbrace  \varphi_q(\xi) ;(1-q^2);  q^2\varphi_q'(\xi) \right\rbrace/\varphi_2(\xi) -\left(12\Phi_3-5\Phi_4 \right) /\varphi_2(\xi)\\
 && =:  7 \Phi_2/\varphi_2(\xi) + \RR '_1 -\RR '_2.
\end{eqnarray}
 According to \eqref{eq1INQ's} and \eqref{Ineq-varphi2} one has 
\begin{align}
\label{eq5INQ's} &\abs{\RR '_1 }\leq 16\cdot 15 \cdot 3!\cdot 4^{2\overline{H}-1}e^{-2\xi}\cdot (1-\delta)^{-\overline{H}} e^{(2-\delta)\xi}  =: M_{3}e^{-\delta\xi} ,\\
\nonumber & \\
\nonumber & \abs{\RR '_2 }\leq \left[ 12\abs{9\varphi_3'(\xi) -\varphi_3(\xi)}+5\abs{16\varphi_4'(\xi) -\varphi_4(\xi)} \right]/\varphi_2(\xi)\\
\label{eq6INQ's}  & \leq \left(108\cdot 3^{\overline{H}-1}+12\cdot 3^{\overline{H}}+80\cdot 4^{\overline{H}-1}+5\cdot 4^{\overline{H}}\right)e^{-2\xi}\cdot (1-\delta)^{-\overline{H}} e^{(2-\delta)\xi} =:M_{4}e^{-\delta\xi} ,
\end{align} 
while 
\begin{align}
 \nonumber & 7 \Phi_2/\varphi_2(\xi) =7 \left(4\frac{\varphi_2'(\xi)}{\varphi_2(\xi)} -1\right)\\
 \nonumber & =7 \left(4\frac{g'(2\xi)}{g'(\xi)}\cdot \frac{g(\xi)}{g(2\xi)} -1\right)=7 \left(4\frac{ (2\xi)^{H-1}K_{H-1}(2\xi)}{\xi^{H-1}K_{H-1}(\xi)}\cdot \frac{\xi^H K_H(\xi)}{(2\xi)^H K_H(2\xi)} -1\right)\\
 \label{eq9INQ's} & =7 \left(2\frac{ K_{H-1}(2\xi)}{ K_H(2\xi)}\cdot \frac{ K_H(\xi)}{  K_{H-1}(\xi)} -1\right)\\ \nonumber & =7 \left(2\frac{ r_{H}(2\xi)}{ r_H(\xi)} -1\right) >7 .
\end{align}
 where $r_H(\xi)=K_{H-1}(\xi)/K_H(\xi)$. In the last  inequality above we used the fact that $r_H(\cdot)$ is strictly increasing on $(0,\infty)$ for $H>1/2$ (see \cite[Lemma 2.2]{gaunt2014inequalities}). 
\paragraph*{\textbf{iv- Conclusions.}}
Putting the previous results \eqref{eq3INQ's},  \eqref{eq8INQ's},  \eqref{eq7INQ's}, \eqref{eq7INQ's-bis}, \eqref{eq5INQ's}, \eqref{eq6INQ's} and \eqref{eq9INQ's}  together with \eqref{DET}, we obtain
\begin{equation*}
\begin{cases}
 \abs{\nabla \varrho_{\theta,3}(\alpha)}=-2\alpha\kappa^3\sigma^4  \xi g(2\xi)g'(\xi) A^*  
 \left[ 
 {\bf\rm{det}} \left\lbrace\CC_1;\CC_3;\CC_5\right\rbrace
 + {\bf\rm{det}} \left\lbrace\CC_2;\CC_3;\CC_5\right\rbrace 
  + {\bf\rm{det}}  \left\lbrace\CC_1;\CC_4;\CC_5\right\rbrace\right]
\\
 ~~\\
 \abs{{\bf\rm{det}} \left\lbrace\CC_{1};\CC_{3};\CC_{5}\right\rbrace}  \leq \frac{M_{1}}{A^*}e^{-(1+\delta)\xi} \mbox{ and } A^*\geq 1, \\
~~\\
 \abs{{\bf\rm{det}} \left\lbrace\CC_2;\CC_{3};\CC_{5}\right\rbrace} \leq M_{2}/(\log(\xi)-M^*),  \\
 ~~\\
{\bf\rm{det}} \left\lbrace\CC_{1};\CC_{4};\CC_{5}\right\rbrace  =   7 \Phi_2/\varphi_2(\xi) + \RR '_1 -\RR '_2 \\
~~\\ 
\abs{\RR'_1}\leq  M_{3}e^{-\delta\xi},~\abs{\RR'_2}\leq  M_{4}e^{-\delta\xi} \mbox{ and } 7 \Phi_2/\varphi_2(\xi) >7,
\end{cases},
\end{equation*}for all $\xi > 1\vee e^{M^*}$. To conclude we choose $(\delta,\xi)\in (0,1)\times (1\vee e^{M^*}, \infty)$  so that
\begin{equation*}
M_1 e^{-(1+\delta)\xi}+\left[ M_{2}/(\log(\xi)-M^*)\right]+(M_3+M_4)e^{-\delta\xi} <6.
\end{equation*}It suffices to choose 
\begin{align}\label{shConst}
\delta\in(0,1)~~\mbox{and}~~\xi=\lambda\alpha > \underset{=:\xi_0}{ \underbrace{1\vee e^{M^*+M_2/3} \vee \frac{1}{\delta}\log \left(\frac{\sum_{j\neq 2}M_j}{3}\right)}}.
\end{align}Under this condition we obtain $\abs{\nabla \varrho_{\theta,3}(\alpha)}<-2\alpha\kappa^3\sigma^4  \xi g(2\xi)g'(\xi) A^* <0$. An equivalent condition for $\alpha$ can be obtained. In fact, we choose  $\alpha >\xi_0/\underline{\lambda}=:\alpha_0$.
 This proves the first part in Lemma \ref{local-inj}, while the second one follows by local univalent theorem (e.g., \cite[p. 2]{parthasarathy1983global}). Therefore, the proof is complete. 

\end{proof}

\subsection*{B2. Evaluating the principal minors of the Jacobian matrix $\nabla \varrho_{\theta,3}(\alpha)$}
In what follows $\langle  K\rangle_\Sigma$ will denote the principal minor of $\nabla \varrho_{\theta,3}(\alpha)$ obtained by omitting rows and columns of indices in $K\subseteq \lbrace 1,2,3\rbrace$. So, the principal minors of order $1\times 1$ are $\langle  2,3\rangle_\Sigma=\Sigma_{1,1}$, $\langle 1,3\rangle_\Sigma=\Sigma_{2,2}$ and  $\langle  1,2\rangle_\Sigma=\Sigma_{3,3}$; while  the principal minors of order $2\times 2$ and order $3\times 3$ are 
\begin{equation*}
\langle 3\rangle_\Sigma=\begin{vmatrix}
\Sigma_{1,1} & \Sigma_{1,2}\\
\Sigma_{2,1} & \Sigma_{2,2}
\end{vmatrix}
~;~~
\langle 2\rangle_\Sigma=\begin{vmatrix}
\Sigma_{1,1} & \Sigma_{1,3}\\
\Sigma_{3,1} & \Sigma_{3,3}
\end{vmatrix}
~;~~
\langle  1\rangle_\Sigma=\begin{vmatrix}
\Sigma_{2,2} & \Sigma_{2,3}\\
\Sigma_{3,2} & \Sigma_{3,3}
\end{vmatrix}
\end{equation*} and
\begin{equation*}
\langle \emptyset\rangle_\Sigma=\begin{vmatrix}
\Sigma_{1,1} & \Sigma_{1,2} & \Sigma_{1,3}\\
\Sigma_{2,1} & \Sigma_{2,2} & \Sigma_{2,3}\\
\Sigma_{3,1} & \Sigma_{3,2} & \Sigma_{3,3}
\end{vmatrix},
\end{equation*}
where $\Sigma_{j,k}$ are the entries of $\nabla \varrho_{\theta,3}(\alpha)$ defined by equations \eqref{eq10}-\eqref{xyzt}.
\begin{lemma}\label{principal-minors}
There exists $\alpha^{*}>0$ depending only on the compact parameter space $$\Theta \subset (0,\infty)\times (1/2,\infty) \times (0,\infty)$$ such that the principal minors of the Jacobian matrix $\nabla \varrho_{\theta,3}(\alpha)$ satisfy
\begin{eqnarray}
\label{princ-minors1} && \langle 2\rangle_\Sigma,~\langle 3\rangle_\Sigma,~\langle  2,3\rangle_\Sigma  > 0~~\mbox{and}~~~\langle 1,2\rangle_\Sigma,~\langle \emptyset\rangle_\Sigma <0,\\
&& \nonumber \\
\label{princ-minors2}  && \begin{cases}
\langle  1\rangle_\Sigma, ~\langle  1,3\rangle_\Sigma     <0 &,~~\mbox{if }\lambda < e^{(\psi(H+1/2)+\psi(H))/2 }\\
\langle  1\rangle_\Sigma,~\langle  1,3\rangle_\Sigma >0 &,~~\mbox{if }\lambda >   e^{(\psi(H+1/2)+\psi(H))/2 }
\end{cases},
\end{eqnarray}  for all $\alpha \in(\alpha^{*},\infty)$.
\end{lemma}
\begin{proof}
First, observe that the unique principal minor of order $3\times 3$ is the determinant of the Jacobian matrix $\nabla \varrho_{\theta,3}(\alpha)$. It is already investigated and shown to be negative when $\alpha >\alpha_0$ (see Lemma \ref{local-inj}). Throughout this proof $\delta \in (2/3,1)$ and the supremum with respect to $H$ and/or $\lambda$ should be understoud as being taken over their own spaces (say $H\in [\underline{H},\overline{H}]$, $\lambda\in [\underline{\lambda},\overline{\lambda}]$ and $(H,\lambda)\in [\underline{H},\overline{H}]\times [\underline{\lambda},\overline{\lambda}]$).
\paragraph*{\textbf{i-Evaluating the principal minors of order $1\times 1$.}} First observe that $\langle  2,3\rangle_\Sigma=\Sigma_{1,1}=2\kappa(g_H(2\xi)-4g_H(\xi)+3g_H(0))$, where  $\kappa>0$ is independent of $\xi=\lambda\alpha$ (see notations \eqref{pre-xyzt15} and \eqref{xyzt}). Letting $\alpha$ to vary in $(0,\infty)$ and consider the function $ \xi\mapsto\langle  2,3\rangle_\Sigma =: f_0(\xi)$. This function is continuously differentiable with derivative
\begin{eqnarray*}
 \partial_\xi f_0(\xi) &=& 2\kappa\left(-2(2\xi)^HK_{H-1}(2\xi)+4\xi^HK_{H-1}(\xi)\right)\\
&=& 8\kappa\xi \left( g_{H-1}(\xi)-g_{H-1}(2\xi)\right)>0
\end{eqnarray*}and $\lim_{\xi\downarrow 0}f_0(\xi)=0$. In the above inequality we used the fact that $\xi\mapsto g_{H-1}(\xi)$ is strictly decreasing on $(0,\infty)$. Thereby, $f_0$ is  strictly increasing and $\langle  2,3\rangle_\Sigma=f_0(\xi)> 0$, $\forall \xi >0$. \\
The principal minor $\langle 1,3\rangle_\Sigma$ can be rewritten as (see notations \eqref{pre-xyzt16} and \eqref{xyzt})
\begin{eqnarray*}
\langle 1,3\rangle_\Sigma &=& \Sigma_{2,2} =\kappa\sigma^2 \lbrace -4 (u_0\bm{\xi}+u^*)+\underset{=:\RR_1(\xi)}{\underbrace{\bm{\xi}(z-4y+7x)+(Cz-4By+7Ax)}}\rbrace\\
&=& \kappa\sigma^2 u_0 \lbrace -4 [-\log(\lambda^2)+\psi(H+1/2)+\psi(H)]+\RR_1(\xi)/u_0\rbrace.
\end{eqnarray*} It is not hard to see that
\begin{eqnarray*}
 \abs{\frac{\RR_1(\xi)}{u_0}} &\leq & \frac{1}{u_0} \left[ \vert \bm{\xi}\vert(z+4y+7x)+(Cz+4By+7Ax)\right]\\
 &\leq &  \frac{x}{u_0} (12 \vert\bm{\xi}\vert +C+4B+7A)\\
  &\leq &  \frac{12g_H(\xi)}{u_0} \left[ \vert -\log(2\lambda^2)+\psi(H+1/2)\vert+  \log(3\xi)+M(6H+1) \right]\\
  &\leq &   12(1-\delta)^{-H}e^{-\delta \xi} \left[ \vert -\log(2\lambda^2)+\psi(H+1/2)\vert+  \log(3\xi)+M(6H+1)\right]\\
  &\leq &  e^{-\delta \xi}(L_0+L_1\log(\xi)),~~\xi=\lambda\alpha,
\end{eqnarray*}for all $\xi>e^M$, where $M$ is the numerical constant given in \eqref{M-const} and used in  {\rm{\textbf{(v)}}},{\rm{\textbf{(vi)}}}-Lemma \ref{LEM0-App}, $L_1=12(1-\delta)^{-\overline{H}}$ and
\begin{eqnarray*}
L_0 &=&\sup_{H,\lambda} \left\lbrace 12(1-\delta)^{-H}\left[ \vert -\log(2\lambda^2)+\psi(H+1/2)\vert+  \log(3)+M(6H+1)\right] \right\rbrace.
\end{eqnarray*}In the first four inequalities we used the following facts: $g_H(q\xi),R_H(q\xi)>0$, for all $\xi>1\vee e^M$ and $q=1,2,3$, the function $\xi\mapsto g_H(q\xi)$ is decreasing, {\rm{\textbf{(vi)}}}-Lemma \ref{LEM0-App} and the first part in \eqref{eq1INQ's}, respectively. We can find some $\xi_1$ (depending only on the parameter space) so that 
\begin{equation*}
\begin{cases}
\langle  1,3\rangle_\Sigma <0 &,~~\mbox{if }\lambda < e^{(\psi(H+1/2)+\psi(H))/2 }\\
\langle  1,3\rangle_\Sigma >0 &,~~\mbox{if }\lambda >   e^{(\psi(H+1/2)+\psi(H))/2 }
\end{cases},~~~~\mbox{for all }\xi>\xi_1.
\end{equation*}Formally, $\xi_1$ can be defined as
\begin{eqnarray*}
\xi_1 &:= &\inf\left\lbrace \xi>e^M \vee e^{1/\delta} : e^{-\delta \xi}(L_0+L_1\log(\xi)) < 4\abs{\psi(H_*+1/2)+\psi(H_*)-\log(\lambda_*^2)} \right\rbrace,\\
& \mbox{with} & (H_*,\lambda_*) := \underset{H,\lambda}{\arg\min}   \abs{\psi(H+1/2)+\psi(H)-\log(\lambda^2)}
\end{eqnarray*} Note that $\xi\mapsto e^{-\delta \xi}(L_0+L_1\log(\xi))$ is strictly decreasing on $(e^{1/\delta},\infty)$.
\\
The principal minor $\langle 1,2\rangle_\Sigma$ has the form (see notations \eqref{pre-xyzt17} and \eqref{xyzt})
\begin{eqnarray*}
\langle 1,2\rangle_\Sigma &=& \Sigma_{3,3} =\kappa\sigma^2  \lbrace \bm{\zeta}u_0 +\underset{=:\RR_2(\xi)}{\underbrace{ \bm{\zeta}(t-4z+6y-4x)-4\alpha(dt-4cz+6by-4ax)}}\rbrace\\
&=& \kappa\sigma^2 u_0 \left( -\frac{2H}{\lambda}+\RR_2(\xi)/u_0\right).
\end{eqnarray*}Let $H>1/2$. Since $\xi\mapsto r_H(\xi)=K_{H-1}(\xi)/K_H(\xi)$ is continuous and increasing on  $(0,\infty)$ (e.g., \cite[p. 376]{gaunt2014inequalities}) with 
\begin{equation*}
\lim_{\xi\downarrow 0}r_H(\xi)=0,~~ \lim_{\xi\uparrow \infty}r_H(\xi)=1
\end{equation*}we get
\begin{equation}\label{r*-const}
 0< r^*:=\inf_H r_H(1)<a<b<c<d<1,~~\mbox{ for all }\xi>1.
\end{equation}
 Moreover,
\begin{eqnarray*}
 \abs{\frac{\RR_2(\xi)}{u_0}} &\leq & \frac{1}{u_0} \left[ \vert \bm{\zeta}\vert (t+4z+6y+4x)+\frac{4}{\lambda}\xi(dt+4cz+6by+4ax)\right]\\
 &\leq &  \frac{30x}{u_0} \left( \frac{H}{\lambda}+\frac{2\xi}{\lambda}\right)\\
  &\leq &   30\lambda^{-1}(1-\delta)^{-H}e^{-\delta \xi} \left( H+ 2\xi\right)\\
  &\leq &  e^{-\delta \xi}(L_2+L_3\xi),
\end{eqnarray*}for all $\xi>1$, where $L_2=\sup_{H,\lambda} \lbrace 30H\lambda^{-1}(1-\delta)^{-H}\rbrace$ and 
$L_3=\sup_{H,\lambda} \lbrace 60\lambda^{-1}(1-\delta)^{-H}\rbrace$. Once again, we  can find some $\xi_2$ (depending only on the parameter space) so that $\langle 1,2\rangle_\Sigma <0$, for all $\xi>\xi_2$.
\paragraph*{\textbf{ii-Evaluating the principal minors of order $2\times 2$.}} 
In the sequel we let $\xi>1\vee e^{M^*}$ so that $A^*\geq 1$ (as shown in \eqref{eq7INQ's}).
Keep in mind notaions \eqref{pre-xyzt15}-\eqref{xyzt} and set $\widetilde{\Sigma}_{1,1}=(y-4x+3u_0)$, $\widetilde{\Sigma}_{2,1}=(z-4y+7x-4u_0)$, $\widetilde{\Sigma}_{3,1}=(t-4z+6y-4x+u_0)$. Then, by straightforward computations we get
\begin{align*}
\langle 3 \rangle_\Sigma  &= 
  \begin{vmatrix}
\displaystyle
2\kappa\widetilde{\Sigma}_{1,1} &~~ 2 \kappa\sigma^2[\bm{\xi} \widetilde{\Sigma}_{1,1}+ (By-4Ax+3u^*)] 
\\
\displaystyle
\kappa\widetilde{\Sigma}_{2,1} &~~  \kappa\sigma^2[\bm{\xi} \widetilde{\Sigma}_{2,1}+(Cz-4By+7Ax-4u^*)] 
\end{vmatrix}\\
&= 2\kappa^2\sigma^2\begin{vmatrix}
\displaystyle
\widetilde{\Sigma}_{1,1} &~~  \bm{\xi} \widetilde{\Sigma}_{1,1}+ (By-4Ax+3u^*)  
\\
\displaystyle
\widetilde{\Sigma}_{2,1} &~~  \bm{\xi} \widetilde{\Sigma}_{2,1}+(Cz-4By+7Ax-4u^*)  
\end{vmatrix}\\
&= 2\kappa^2\sigma^2\begin{vmatrix}
\displaystyle
\widetilde{\Sigma}_{1,1} &~~    By-4Ax+3u^*  
\\
\displaystyle
\widetilde{\Sigma}_{2,1} &~~   Cz-4By+7Ax-4u^* 
\end{vmatrix}\\
&= 2\kappa^2\sigma^2\begin{vmatrix}
\displaystyle
y-4x+3u_0 &~~    By-4Ax+3u^* 
\\
\displaystyle
z-4y+7x-4u_0 &~~   Cz-4By+7Ax-4u^* 
\end{vmatrix}\\
&= 2\kappa^2\sigma^2\begin{vmatrix}
\displaystyle
y-4x+3u_0 &~~    (B-A)y-3A^* 
\\
\displaystyle
z-4y+7x-4u_0 &~~   (C-A)z-4(B-A)y+4A^* 
\end{vmatrix},~~~A^*=Au_0-u^*\\
&= 2\kappa^2\sigma^2\begin{vmatrix}
\displaystyle
(y-4x)+3u_0 &~~    (B-A)y-3A^* 
\\
\displaystyle
(z-9x)+8u_0 &~~   (C-A)z-8A^*
\end{vmatrix}\\
&= 2\kappa^2\sigma^2 \cdot \left\lbrace 
\begin{vmatrix}
\displaystyle
(y-4x)  &~~    (B-A)y  
\\
\displaystyle
(z-9x) &~~   (C-A)z
\end{vmatrix} 
+\begin{vmatrix}
\displaystyle
(y-4x)  &~~     -3A^* 
\\
\displaystyle
(z-9x)  &~~    -8A^*
\end{vmatrix} 
+\begin{vmatrix}
\displaystyle
 3u_0 &~~    (B-A)y  
\\
\displaystyle
 8u_0 &~~   (C-A)z 
\end{vmatrix} 
\right\rbrace 
\end{align*} 
\begin{align*} 
&= 2\kappa^2\sigma^2 A^*x \cdot \left\lbrace 
\begin{vmatrix}
\displaystyle
(y/x-4)  &~~    (B-A)y/A^* 
\\
\displaystyle
(z/x-9) &~~   (C-A)z/A^* 
\end{vmatrix} 
+\begin{vmatrix}
\displaystyle
(y/x-4)  &~~     -3  
\\
\displaystyle
(z/x-9)  &~~    -8 
\end{vmatrix} 
+\begin{vmatrix}
\displaystyle
 3u_0/A^* &~~    (B-A)y/x 
\\
\displaystyle
 8u_0/A^* &~~   (C-A)z/x
\end{vmatrix} 
\right\rbrace\\
&=: 2\kappa^2\sigma^2 A^*x \sum_{j=3}^5 \RR_j(\xi).
\end{align*} We can evaluate the last three deteminants $\RR_j(\xi)$, in terms of $g(\xi):=g_H(\xi)$ and $\varphi_q(\xi)=g(q\xi)/g(\xi)$, $R(\xi):=R_H(q\xi)$, $q=1,2,3$ as 
\begin{eqnarray*}
\vert  \RR_3(\xi)\vert &=& \abs{\frac{z}{A^*}(C-A)(y/x-4)  -  \frac{y}{A^*}(B-A) (z/x-9)}\\
&\leq & \abs{\frac{g(3\xi)}{A^*}[R(3\xi)-R(\xi)]}\left(\varphi_2(\xi)+4\right)+\abs{\frac{g(2\xi)}{A^*}[R(2\xi)-R(\xi)]}\left(\varphi_3(\xi)+9\right)\\
&\leq & \frac{M_0}{A^*}g(2\xi)\left[ \varphi_2(\xi)+4 + \varphi_3(\xi)+9\right],~~\mbox{(By \eqref{eq2INQ's})}\\
&\leq & \frac{M_0}{A^*} u_0(1-\delta)^{-H}e^{-2\delta\xi} (13+2^H e^{-\xi}+3^H e^{-2\xi} ),~~\mbox{(By \eqref{eq1INQ's})}\\
&\leq &  M_0 e^{-2\delta\xi} \sup_H\left\lbrace u_0(1-\delta)^{-H} (13+2^{H}  +3^{H} )\right\rbrace ,~~~\mbox{for all }\xi>e^{M^*}.
\end{eqnarray*} On the other hand, $\RR_4(\xi)  = -8(y/x-4) +3(z/x-9)  =5-8\varphi_2(\xi)+3\varphi_3(\xi)$ with
\begin{eqnarray*}
  \abs{-8\varphi_2(\xi)+3\varphi_3(\xi)}  &\leq & 8\varphi_2(\xi)+3\varphi_3(\xi)\\
&\leq &  2^{H+3} e^{-\xi}+3^{H+1} e^{-2\xi},~~\mbox{(By the second part in \eqref{eq1INQ's})}\\
&\leq & ( 2^{\overline{H}+3}+3^{\overline{H}+1}) e^{-\xi}, ,~~~\mbox{for all }\xi>1.
\end{eqnarray*}and for $\xi>1\vee e^{M^*}$ (see \eqref{eq7INQ's})  we have
\begin{eqnarray*}
\vert  \RR_5(\xi)\vert &=& \frac{u_0}{A^*} \abs{3(C-A)\frac{z}{x}-8(B-A)\frac{y}{x} }\\
&\leq & u_0 [3\varphi_3(\xi)\abs{R(3\xi)-R(\xi)}+8\varphi_2(\xi)\abs{R(2\xi)-R(\xi)} ]\\
&\leq & u_0 [\log(3) +3M(2H+1)]\cdot \left[ 8\varphi_2(\xi)+3\varphi_3(\xi)\right],~~\mbox{(By {\rm{\textbf{(iv)}}}-Lemma \ref{LEM0-App})}\\
&\leq & u_0 [\log(3) +3M(2H+1)]\cdot \left[ 2^{H+3} e^{-\xi}+3^{H+1} e^{-2\xi}\right]\\
&\leq & L_4e^{-\xi},
\end{eqnarray*}where $L_4=\sup_H \left\lbrace u_0 [\log(3) +3M(2H+1)]\cdot \left[ 2^{H+3}  +3^{H+1}  \right] \right\rbrace$. Putting the previous estimates of $\RR_j(\xi)$ together, we can find some $\xi_3$  (depending only on the parameter space) so that $\langle 3\rangle_\Sigma >0$, for all $\xi>\xi_3$.\bigskip

Similarly, we evaluate the principal minor $\langle 2\rangle_\Sigma$ as
 \begin{align*}
\langle 2 \rangle_\Sigma  &= 
  \begin{vmatrix}
\displaystyle
2\kappa\widetilde{\Sigma}_{1,1}  &~~ 2 \kappa\sigma^2[\bm{\zeta} \widetilde{\Sigma}_{1,1} -2\alpha(by-2ax)]
\\
\displaystyle
\kappa\widetilde{\Sigma}_{3,1} &~~  \kappa\sigma^2[\bm{\zeta} \widetilde{\Sigma}_{3,1}-4\alpha(dt-3cz+3by-ax)]
\end{vmatrix}\\
&= 2 \kappa^2\sigma^2
  \begin{vmatrix}
\displaystyle
 \widetilde{\Sigma}_{1,1}  &~~  \bm{\zeta} \widetilde{\Sigma}_{1,1} -2\alpha(by-2ax) 
\\
\displaystyle
\widetilde{\Sigma}_{3,1} &~~    \bm{\zeta} \widetilde{\Sigma}_{3,1}-4\alpha(dt-3cz+3by-ax) 
\end{vmatrix}\\
&= 2 \kappa^2\sigma^2
  \begin{vmatrix}
\displaystyle
 \widetilde{\Sigma}_{1,1}  &~~    -2\alpha(by-2ax) 
\\
\displaystyle
\widetilde{\Sigma}_{3,1} &~~     -4\alpha(dt-3cz+3by-ax) 
\end{vmatrix}\\
\end{align*}
\begin{align*}
&= -4\alpha \kappa^2\sigma^2
  \begin{vmatrix}
\displaystyle
  y-4x+3u_0 &~~     by-2ax
\\
\displaystyle
t-4z+6y-4x+u_0 &~~     2dt-6cz+6by-2ax
\end{vmatrix}\\
&= -4\alpha \kappa^2\sigma^2 \left\lbrace
  \begin{vmatrix}
  y-4x  &~~     by 
\\
t-4z+6y-4x  &~~     2dt-6cz+6by 
\end{vmatrix}
+  \begin{vmatrix}
  y-4x  &~~      -2ax
\\
t-4z+6y-4x  &~~      -2ax
\end{vmatrix}\right.\\
& ~~~~~~~~~~~~~~~~~~~~~~~~~~~~\left. +  \begin{vmatrix}
   3u_0 &~~     by 
\\
 u_0 &~~     2dt-6cz+6by 
\end{vmatrix}
+  \begin{vmatrix}
  3u_0 &~~      -2ax
\\
 u_0 &~~     -2ax
\end{vmatrix}
\right\rbrace
\\ 
&= -4\alpha \kappa^2\sigma^2 \cdot (au_0x) \left\lbrace
 \frac{1}{au_0} \begin{vmatrix}
  y/x-4  &~~     by/x 
\\
t-4z+6y-4x  &~~     2dt-6cz+6by 
\end{vmatrix}
- \frac{2}{u_0} \begin{vmatrix}
  y-4x  &~~      1
\\
t-4z+6y-4x  &~~      1
\end{vmatrix}\right.\\
& ~~~~~~~~~~~~~~~~~~~~~~~~~~~~\left. + \frac{1}{a} \begin{vmatrix}
   3  &~~     by/x
\\
 1 &~~     (2dt-6cz+6by)/x 
\end{vmatrix}
-2  \begin{vmatrix}
  3  &~~     1 
\\
 1 &~~     1 
\end{vmatrix}
\right\rbrace
\\
&=: -4 au_0\alpha \kappa^2\sigma^2 x\left( \sum_{j=6}^8\RR_j(\xi) -4 \right),~~~\xi=\lambda\alpha.
\end{align*}The terms $\RR_j(\xi)$ can be controlled to be small. In fact, by using \eqref{r*-const} and the fact: $t<z<y<x$ (inherited from the monotony of $\xi\mapsto g_H(\xi)$) we obtain
\begin{eqnarray*}
\vert \RR_6(\xi) \vert &= & \frac{1}{r^* u_0}\abs{(2d-b)\frac{yt}{x}+(-6c+4b)\frac{yz}{x}-8dt+24cz-20by}\\
  &\leq  & \frac{1}{r^* u_0}\left(  \frac{y}{x}(2t+6z)+8t+24z+20y\right)\\
    &\leq  & \frac{1}{r^* u_0}\left(  8y+8t+24z+20y\right) \leq  \frac{60y}{r^* u_0}\\
&\leq &    \frac{60}{r^* u_0}g(2\xi) \leq \frac{60}{r^*}(1-\delta)^{-\overline{H}}e^{-2\delta\xi} ,~~\mbox{(By the first part in \eqref{eq1INQ's})},
\end{eqnarray*}
\begin{eqnarray*}
\vert \RR_7(\xi) \vert &= & \frac{2}{ u_0}\abs{(y-4x)-(t-4z+6y-4x)}=\frac{2}{ u_0}\abs{-t+4z-5y}\\
&\leq & \frac{20 g(2\xi)}{ u_0}\leq  20(1-\delta)^{-\overline{H}}e^{-2\delta\xi} ,~~\mbox{(By the first part in \eqref{eq1INQ's})}
\end{eqnarray*}and 
\begin{eqnarray*}
\vert \RR_8(\xi) \vert &= & \frac{1}{a}\abs{6 \frac{dt}{x} -18\frac{cz}{x}+17\frac{by}{x}} \leq \frac{41y}{ax}\\
&\leq & \frac{41 }{ r^*} \varphi_2(\xi)\leq \frac{41 }{ r^*} 2^{\overline{H}}e^{-\xi}
 ,~~\mbox{(By the second part in \eqref{eq1INQ's})}
\end{eqnarray*} Combining the previous estimates of $\RR_j(\xi)$ we can find some $\xi_4>0$ (depending only on the parameter space), so that $\langle 2\rangle_\Sigma >0$, for all $\xi>\xi_4$.\bigskip

Now, let us examine the principal minor $\langle 1\rangle_\Sigma$. Recall that  $\widetilde{\Sigma}_{2,1}=(z-4y+7x-4u_0)$, $\widetilde{\Sigma}_{3,1}=(t-4z+6y-4x+u_0)$, $\bm{\xi}=-\log \left(2\lambda^2\right)+\psi(H+1/2)$, $\bm{\zeta}=-2H/\lambda$ and write 
\begin{eqnarray*}
&& \bm{\xi} \widetilde{\Sigma}_{2,1}+(Cz-4By+7Ax-4u^*) =-4(u_0\bm{\xi} +u^*)+\underset{=:f_1(\xi)}{\underbrace{\bm{\xi} (z-4y+7x)+(Cz-4By+7Ax)}} \\
&& \bm{\zeta} \widetilde{\Sigma}_{2,1}-\alpha(3cz-8by+7ax)  =-4u_0\bm{\zeta}-7a\alpha x+\underset{=:f_2(\xi)}{\underbrace{\bm{\zeta}(z-4y+7x)  -\alpha(3cz-8by)}}\\
&& ~~\\
&& \bm{\xi} \widetilde{\Sigma}_{3,1}+(Dt-4Cz+6By-4Ax+u^*)= (u_0\bm{\xi}+u^*) \\
&& ~~~~~~~~~~~~~~~~~~~~~~~~~~~~~~~~~~~~~~ +\underset{=:f_3(\xi)}{\underbrace{\bm{\xi}(t-4z+6y-4x)+(Dt-4Cz+6By-4Ax) }} \\
&& \bm{\zeta} \widetilde{\Sigma}_{3,1}-4\alpha(dt-3cz+3by-ax) = u_0\bm{\zeta} +4a\alpha x+\underset{=:f_4(\xi)}{\underbrace{\bm{\zeta}(t-4z+6y-4x) -4\alpha(dt-3cz+3by)}}
\end{eqnarray*} Then the principal minor $\langle 1\rangle_\Sigma$ simplifies to
\begin{align*}
\langle 1 \rangle_\Sigma  &= 
  \kappa^2\sigma^4 \begin{vmatrix}
\displaystyle
    -4( u_0\bm{\xi}+u^*)+ f_1(\xi) &~~   -4u_0\bm{\zeta}-7a\alpha x+f_2(\xi)
\\
\displaystyle
    (u_0\bm{\xi}+u^*)  + f_3(\xi) &~~   u_0\bm{\zeta} +4a\alpha x+f_4(\xi)
\end{vmatrix}\\
 &=\kappa^2\sigma^4\left\lbrace 
 \begin{vmatrix}
    -4( u_0\bm{\xi}+u^*)  &~~   -4u_0\bm{\zeta}-7a\alpha x 
\\
    (u_0\bm{\xi}+u^*)    &~~   u_0\bm{\zeta} +4a\alpha x 
\end{vmatrix} 
+\begin{vmatrix}
    -4(u_0\bm{\xi} +u^*)  &~~    f_2(\xi)
\\
    (u_0\bm{\xi}+u^*)    &~~    f_4(\xi)
\end{vmatrix}\right.\\
& \left. ~~~~~~~~~~~~~~~~~~~~~~~~~~~~~~~~~~+\begin{vmatrix}
     f_1(\xi) &~~   -4u_0\bm{\zeta}-7a\alpha x 
\\
      f_3(\xi) &~~   u_0\bm{\zeta} +4a\alpha x 
\end{vmatrix}
+\begin{vmatrix}
     f_1(\xi) &~~    f_2(\xi)
\\
      f_3(\xi) &~~   f_4(\xi)
\end{vmatrix}
\right\rbrace\\
 &=\kappa^2\sigma^4\cdot (a\alpha x)\left\lbrace 
 \begin{vmatrix}
    -4( u_0\bm{\xi}+u^*)  &~~    -7 
\\
    (u_0\bm{\xi}+u^*)    &~~   4 
\end{vmatrix} 
+\begin{vmatrix}
    -4(u_0\bm{\xi} +u^*)  &~~    f_2(\xi)/(a\alpha x)
\\
    (u_0\bm{\xi}+u^*)    &~~    f_4(\xi)/(a\alpha x)
\end{vmatrix}\right.\\
& \left. ~~~~~~~~~~~~~~~~~~~~~~~~~~~~~~~~~~+\begin{vmatrix}
     f_1(\xi)/x &~~   -4u_0\bm{\zeta}/(a\alpha )-7x 
\\
      f_3(\xi)/x &~~   u_0\bm{\zeta}/(a\alpha ) +4x
\end{vmatrix}
+\begin{vmatrix}
     f_1(\xi) &~~    f_2(\xi)/(a\alpha x)
\\
      f_3(\xi) &~~   f_4(\xi)/(a\alpha x)
\end{vmatrix}
\right\rbrace\\
&=a\alpha\kappa^2\sigma^4  x \left\lbrace 
 -9 (u_0\bm{\xi} +u^*)
+(u_0\bm{\xi} +u^*)\begin{vmatrix}
    -4  &~~    f_2(\xi)/(a\alpha x)
\\
    1   &~~    f_4(\xi)/(a\alpha x)
\end{vmatrix}\right.\\
& \left. ~~~~~~~~~~~~~~~~~~~~~~~~~~~~~~~~~~+\begin{vmatrix}
     f_1(\xi)/x &~~   -4u_0\bm{\zeta}/(a\alpha )-7 x
\\
      f_3(\xi)/x &~~   u_0\bm{\zeta}/(a\alpha ) +4x
\end{vmatrix}
+\begin{vmatrix}
     f_1(\xi) &~~    f_2(\xi)/(a\alpha x)
\\
      f_3(\xi) &~~   f_4(\xi)/(a\alpha x)
\end{vmatrix}
\right\rbrace\\
&=:\kappa^2\sigma^4 \left(-9 (u_0\bm{\xi} +u^*)+\sum_{j=9}^{11}\RR_j(\xi) \right).
\end{align*}Simple computations yield
\begin{eqnarray*}
\vert \RR_9(\xi)\vert &\leq & \vert u_0\bm{\xi} +u^*\vert \cdot \frac{\vert 4f_4(\xi)+f_2(\xi) \vert }{a\alpha x}\\
&\leq & \frac{\vert u_0\bm{\xi} +u^*\vert}{r^*}\left( \frac{\vert \bm{\zeta}(4t-15z+20y-9x)\vert }{\alpha x}+\frac{\vert \alpha(16dt-45cz+40by)\vert}{\alpha x}\right)\\
&\leq & \frac{\vert u_0\bm{\xi} +u^*\vert}{r^*}\left[ \frac{\lambda\vert \bm{\zeta}\vert}{\xi}\left(\frac{4t}{x}+\frac{15z}{x}+\frac{20y}{x}+9\right)+  \left(\frac{16dt}{x}+\frac{45cz}{x}+\frac{40by}{x}\right) \right],~~\xi=\lambda\alpha\\
&\leq & \frac{\vert u_0\bm{\xi} +u^*\vert}{r^*}\left( 48\frac{\lambda\vert \bm{\zeta}\vert}{\xi }+  101\frac{y}{x}  \right)=\frac{\vert u_0\bm{\xi} +u^*\vert}{r^*}\left( 48\frac{\lambda\vert \bm{\zeta}\vert}{\xi }+  101\varphi_2(\xi)  \right)\\
&\leq & \frac{\vert u_0\bm{\xi} +u^*\vert}{r^*}\left( 48\frac{\lambda\vert \bm{\zeta}\vert}{\xi }+  101\cdot 2^H e^{-\xi}  \right) \leq   \frac{L_5}{\xi},
\end{eqnarray*}where $L_5=\sup_{H,\lambda} \left\lbrace \frac{\vert u_0\bm{\xi} +u^*\vert}{r^*}\left( 48 \lambda\vert \bm{\zeta}\vert +  101\cdot 2^H    \right) \right\rbrace$. Also we have
\begin{eqnarray*}
\vert \RR_{10}(\xi)\vert & = & \abs{\frac{f_1(\xi)}{x} \left(\frac{u_0\bm{\zeta}}{a\alpha } +4x\right)    +  \frac{f_3(\xi)}{x}\left(4\frac{u_0\bm{\zeta}}{a\alpha }+7 x\right)}\\
& = & \frac{u_0\lambda \vert\bm{\zeta}\vert}{a\xi }\left(\frac{\vert f_1(\xi)\vert}{x}+4\frac{\vert f_3(\xi)\vert}{x}\right)+4 \vert f_1(\xi)\vert +7 \vert f_3(\xi)\vert,~~\xi=\lambda\alpha
\end{eqnarray*}with
\begin{eqnarray*} 
\frac{\vert f_1(\xi)\vert}{x} &\leq & \vert \bm{\xi}\vert \left( \frac{z}{x}+4\frac{y}{x}+7\right)+\left( \frac{\vert C\vert z}{x}+4\frac{\vert B\vert y}{x}+7\vert A\vert \right)\\
&\leq & 12 \left[\vert \bm{\xi}\vert +\log(4\xi)+M(6H+1) \right],\\
\frac{\vert f_3(\xi)\vert}{x} &\leq & \vert \bm{\xi}\vert \left( \frac{t}{x}+4\frac{z}{x}+6\frac{y}{x}+4\right)+\left( \frac{\vert D\vert t}{x}+4\frac{ \vert C\vert z}{x}+6\frac{\vert B\vert y}{x}+4\vert A\vert \right)\\
&\leq & 15 \left[\vert \bm{\xi}\vert +\log(4\xi)+M(6H+1) \right],\\
&&~~\\
\vert f_1(\xi)\vert & = & \frac{\vert f_1(\xi)\vert}{x} \cdot x \leq   12 u_0\left[\vert \bm{\xi}\vert +\log(4\xi)+M(6H+1) \right](1-\delta)^{-H}e^{-\delta \xi}\\
\vert f_3(\xi)\vert &\leq & 15 u_0\left[\vert \bm{\xi}\vert +\log(4\xi)+M(6H+1) \right](1-\delta)^{-H}e^{-\delta \xi}
\end{eqnarray*}It follows that 
\begin{eqnarray*}
\vert \RR_{10}(\xi)\vert & \leq  &  \left[\vert \bm{\xi}\vert +\log(4\xi)+M(6H+1) \right]\left(
72\frac{u_0\lambda \vert\bm{\zeta}\vert}{r^* \xi} +153 u_0(1-\delta)^{-H}e^{-\delta \xi} \right)\\
&\leq &  \xi^{-1/3}\left[\vert \bm{\xi}\vert +\log(4\xi)+M(6H+1) \right] \left(
\frac{72u_0\lambda \vert \bm{\zeta}\vert}{r^* } +\frac{153 u_0}{\delta (1-\delta)^{H} } \right)\xi^{-2/3}\\
&\leq &   \left[\vert \bm{\xi}\vert +3\sqrt[3]{4}+M(6H+1) \right] \left(
\frac{72u_0\lambda \vert \bm{\zeta}\vert}{r^* } +\frac{153 u_0}{\delta (1-\delta)^{H} } \right)\xi^{-2/3}\\
&\leq &  L_6 \xi^{-2/3},~~~\mbox{ for all }\xi>1,
\end{eqnarray*}where 
$L_6=\sup_{H,\lambda} \left\lbrace \left[\vert \bm{\xi}\vert +3\sqrt[3]{4}+M(6H+1) \right]\left(
\frac{72u_0\lambda \vert \bm{\zeta}\vert}{r^* } +\frac{153 u_0}{\delta (1-\delta)^{H} } \right)\right\rbrace$.
 Finally, we evaluate the last determinant as
\begin{eqnarray*}
\vert\RR_{11}(\xi)\vert &=& \abs{f_1(\xi) \frac{f_4(\xi)}{a\alpha x}-f_3(\xi)\frac{f_2(\xi)}{a\alpha x}}\\
&\leq &   \frac{\vert f_1(\xi) \vert }{r^*} \left[ \frac{\vert \bm{\zeta}\vert}{\alpha} \left(  \frac{t}{x}+4\frac{z}{x}+6\frac{y}{x}+4\right)+4\left( \frac{dt}{x}+3\frac{cz}{x}+3\frac{b y}{x}\right) \right]  \\
&& ~~~~~~+  \frac{\vert f_3(\xi) \vert }{r^*} \left[ \frac{\vert \bm{\zeta}\vert}{\alpha} \left(  \frac{z}{x}+4\frac{y}{x}+7\right)+ \left( 3\frac{cz}{x}+8\frac{b y}{x}\right) \right]\\
&\leq &   \frac{\vert f_1(\xi) \vert }{r^*} \left( 15\frac{\lambda\vert \bm{\zeta}\vert}{\xi} +28\varphi_2(\xi) \right)   +  \frac{\vert f_3(\xi) \vert }{r^*} \left( 12\frac{\lambda\vert \bm{\zeta}\vert}{\xi} +11\varphi_2(\xi)\right),~~~\varphi_2(\xi)=y/x \\
&\leq &   \frac{\vert f_1(\xi) \vert }{r^*} \left( 15\frac{\lambda\vert \bm{\zeta}\vert}{\xi} +28\cdot 2^H e^{-\xi} \right)   +  \frac{\vert f_3(\xi) \vert }{r^*} \left( 12\frac{\lambda\vert \bm{\zeta}\vert}{\xi} +11\cdot 2^H e^{-\xi}\right)
\end{eqnarray*}
By using the fact that $\log(w)\leq 3w^{1/3}$, for all $w>0$ we obtain
\begin{eqnarray*}
 \vert f_1(\xi)\vert \vee \vert f_3(\xi)\vert &\leq &  15 u_0\left[\vert \bm{\xi}\vert +\log(4\xi)+M(6H+1) \right](1-\delta)^{-H}e^{-\delta \xi}\\
&\leq & 15 u_0\left[\vert \bm{\xi}\vert +3\sqrt[3]{4}+M(6H+1) \right](1-\delta)^{-H} \xi^{1/3}e^{-\delta \xi}\\
&\leq & L_7\xi^{1/3}e^{-\delta \xi},~~~\mbox{for all }\xi>1,
\end{eqnarray*} where $$L_7:=\sup_{H,\lambda} \left\lbrace 15 u_0\left[\vert \bm{\xi}\vert +3\sqrt[3]{4}+M(6H+1) \right](1-\delta)^{-H} \right\rbrace.$$
 Hence,  
\begin{eqnarray*}
\vert\RR_{11}(\xi)\vert &\leq & \frac{L_7}{r^*} \left(27\frac{\lambda\vert \bm{\zeta}\vert}{\xi} +39\cdot 2^H e^{-\xi} \right)\xi^{1/3}e^{-\delta \xi}\\
&& \leq  L_8\xi^{-2/3}e^{-\delta \xi},
\end{eqnarray*}
where $L_8:=\sup_{H,\lambda} \left\lbrace \frac{L_7}{r^*} \left(27 \lambda\vert \bm{\zeta}\vert  +39\cdot 2^H     \right) \right\rbrace$. According to the previous estimates of $\RR_j(\xi)$, $j=9,10,11$ we can find some $\xi_5>0$ (depending only on the
parameter space), so that for all $\xi>\xi_5$ we have
\begin{equation*}
\begin{cases}
\langle  1\rangle_\Sigma <0 &,~~\mbox{if }u_0\bm{\xi} +u^*>0  \\
\langle  1\rangle_\Sigma >0 &,~~\mbox{if }u_0\bm{\xi} +u^*<0
\end{cases},
\end{equation*}Or equivalently,
\begin{equation*}
\begin{cases}
\langle  1\rangle_\Sigma <0 &,~~\mbox{if }\lambda < e^{(\psi(H+1/2)+\psi(H))/2 }\\
\langle  1\rangle_\Sigma >0 &,~~\mbox{if }\lambda >   e^{(\psi(H+1/2)+\psi(H))/2 }
\end{cases}
\end{equation*}
\paragraph*{\textbf{iii-Conclusion.}}  All the proper principal minors are evaluated according to the values of $\xi =\lambda\alpha $ in the intervals $(\xi_j,\infty)$, $j=1,2,3,4,5$. It is then clear that both \eqref{princ-minors1} and \eqref{princ-minors2} hold for all $\alpha\in(\alpha^*,\infty)$, where $\alpha^{*}:=\alpha_0\vee \sup_{1\leq j\leq 5} \lbrace  \xi_j/\underline{\lambda} \rbrace$ and $\alpha_0$ is due to Lemma \ref{local-inj}. The proof is then complete.

\end{proof}
\subsection*{B3. Injectivity of the mapping $\theta \mapsto \B\varrho_{\theta}(\alpha)$}
Let the assumptions of Theorem \ref{TH2} be fulfilled, and let $\theta,\theta'\in \Theta$ such that $\B \varrho_\theta(\alpha)=\B \varrho_{\theta'}(\alpha)$, where $\varrho_\theta(\alpha)=(\gamma_{\theta}(0),\gamma_{\theta}(\alpha),\gamma_{\theta}(2\alpha),\gamma_{\theta}(3\alpha))'$. Set 
\begin{equation}\label{diff-eta_j}
  \gamma_{\theta}(j\alpha)-\gamma_{\theta'}(j\alpha)=\eta_j,~~\mbox{for all }j=0,1,2,3.
\end{equation}Suppose that $\eta_j$   are not all equal, that is, there is some $j^* \neq 0$ such that $\eta_0\neq \eta_{j^*}$. Without loss of generality, say $\eta_0 \neq\eta_1$. The other cases ($j^*=2$ or $3$) can be treated in a similar way to produce the same result. Let $u$ and $v$ be two reals such that
\begin{equation*}
\begin{cases}
u\neq 0, ~v=-u \eta_0/\eta_1,~~\mbox{ if }\eta_1\neq 0 \\[1mm]
v\neq 0 , ~u=0,~~\mbox{ if }\eta_1= 0
\end{cases}
\end{equation*}  It follows that $u(\gamma_{\theta}(0)-\gamma_{\theta'}(0))+v(\gamma_{\theta}(\alpha)-\gamma_{\theta'}(\alpha))=u\eta_0+v\eta_1=0$, or equivalently,
\begin{equation*}
u\gamma_{\theta}(0)+v\gamma_{\theta}(\alpha)=u\gamma_{\theta'}(0)+v\gamma_{\theta'}(\alpha).
\end{equation*}Combining this with the equality $\B \varrho_\theta(\alpha)=\B \varrho_{\theta'}(\alpha)$ yields
\begin{equation}\label{QUQ}
\widetilde{\B}\varrho_\theta(\alpha) =\widetilde{\B} \varrho_{\theta'}(\alpha),~~\widetilde{\B}:=\begin{bmatrix}
\B \\ w
\end{bmatrix},
\end{equation}where $w$ is the row vector defined as $w=(u, v,0 ,0)$. Clearly, the determinant of $\widetilde{\B}$ has the form $\vert\widetilde{\B}\vert=32(u+v)\neq 0$. Thus $\widetilde{\B}$ is invertible and \eqref{QUQ} becomes $\varrho_\theta(\alpha) =\varrho_{\theta'}(\alpha)$. Finally, by injectivity of the mapping $\theta \mapsto\varrho_\theta(\alpha)$ (inherited from the injectivity of $\theta \mapsto\varrho_{\theta,3}(\alpha)$) we get $\theta=\theta '$. We conclude that  $\theta \mapsto \B\varrho_\theta(\alpha)$ is injective on the set of parameters whose covariance vectors cannot differ by a common constant, that is, $\eta_j$ are not all equal. The case where $\eta_j$ are all equal is not possible for large values of $\alpha$, unless all $\eta_j$ are identically zero which in turn implies $\varrho_{\theta,3}(\alpha)=\varrho_{\theta',3}(\alpha)$. Hence, $\theta=\theta'$. This is due to the asymptotic behaviour of $g_H(\cdot)$ at infinity.
  As previously done for $\nabla \varrho_{\theta,3}(\alpha)$, a detailed investigation of the principal minors of $\nabla \B\varrho_\theta(\alpha)=\B\nabla\varrho_\theta(\alpha)$ is possible, and an appropriate use of the existing results on univalence \cite{gale1965jacobian,parthasarathy1983global,olech1991almost,inada1971production} will lead to global injectivity on the parameter space $\Theta$. The next lemma is very crucial for the statement of Theorem \ref{TH3}, since it justifies the existence of the matrix $\CC(\theta)$.
\begin{lemma}\label{LEM-B3}
There is some $\widetilde{\alpha}>0$ (depending only on the parameter space $\Theta$), such that $\abs{\B \nabla \varrho_{\theta}(\alpha)}<0$, for all $\alpha>\widetilde{\alpha}$. Moreover, the mapping $\theta \mapsto \B\varrho_{\theta}(\alpha)$ is locally injective.
\end{lemma}
\begin{proof}\label{proof-LEM-B3}
We are concerned with the Jacobian matrix of $\varrho_{\theta}(\alpha) $ with respect to $\theta$ given as
\begin{equation*}
\displaystyle
\nabla \varrho_{\theta}(\alpha) = \begin{pmatrix}
\displaystyle
\Sigma_{1,1}~~ & \Sigma_{1,2}~~ & \Sigma_{1,3}
\\
\displaystyle
\Sigma_{2,1}~~ & \Sigma_{2,2}~~ & \Sigma_{2,3}
\\
\displaystyle
\Sigma_{3,1}~~ & \Sigma_{3,2}~~ & \Sigma_{3,3} 
\\
\displaystyle
\Sigma_{4,1}~~ & \Sigma_{4,2}~~ & \Sigma_{4,3} \\
\end{pmatrix}
\end{equation*}where its components $\Sigma_{i,j}$ are given in \eqref{pre-xyzt15}-\eqref{pre-xyzt17}, \eqref{xyzt} and 
\begin{equation*}
\begin{cases} 
     \Sigma_{4,1}:=\partial_{\sigma^2} \gamma_\theta(3\alpha)  &= \kappa \left( u - 4t + 6z - 4y + x \right) \\  
    \Sigma_{4,2}:=\partial_H \gamma_\theta(3\alpha)      &= \sigma^2 \kappa \left[ \bm{\xi} (u - 4t + 6z - 4y + x) + \left( Fu - 4Dt + 6Cz - 4By + Ax \right) \right] \\ 
         \Sigma_{4,3}:=\partial_\lambda \gamma_\theta(3\alpha)  &= \sigma^2 \kappa \left[ \bm{\zeta} (u - 4t + 6z - 4y + x) - \alpha \left( 5fu - 16dt + 18cz - 8by + ax \right) \right]  
\end{cases},
\end{equation*}with $u=g_H(5\lambda\alpha)$, $F = R_H(5\lambda\alpha)$ and $f = r_H(5\lambda\alpha)$. Set $\mathbf{P}=\B\nabla\varrho_{\theta}(\alpha)$. It is $3\times 3$-matrix with entries $\mathbf{P}_{i,j}$ given as
\begin{equation*}
\begin{cases}
& \mathbf{P}_{1,1}  = \kappa \left( -2t + 8z - 8y - 8x + 10u_0 \right) \\
&\mathbf{P}_{2,1} = \kappa \left( 2t - 16z + 56y - 112x + 70u_0 \right) \\
&\mathbf{P}_{3,1} = \kappa \left( 2u - 12t + 26z - 16y - 28x + 28u_0 \right) \\[3mm]
&\mathbf{P}_{1,2}  = \sigma^2 \bm{\xi} P_{1,1} + 2\kappa\sigma^2 \left( -Dt + 4Cz - 4By - 4Ax + 5u^* \right) \\
&\mathbf{P}_{2,2} = \sigma^2 \bm{\xi} P_{2,1} + 2\kappa\sigma^2 \left( Dt - 8Cz + 28By - 56Ax + 35u^* \right) \\
&\mathbf{P}_{3,2} = \sigma^2 \bm{\xi} P_{3,1} + 2\kappa\sigma^2 \left( Fu - 6Dt + 13Cz - 8By - 14Ax + 14u^* \right) \\[3mm]
&\mathbf{P}_{1,3}  = \sigma^2 \bm{\zeta} P_{1,1} + 8\kappa\sigma^2\alpha \left( dt - 3cz + 2by + ax \right) \\
&\mathbf{P}_{2,3} = \sigma^2 \bm{\zeta} P_{2,1} - 8\kappa\sigma^2\alpha \left( dt - 6cz + 14by - 14ax \right) \\ 
&\mathbf{P}_{3,3} = \sigma^2 \bm{\zeta} P_{3,1} - 2\kappa\sigma^2\alpha \left( 5fu - 24dt + 39cz - 16by - 14ax \right)
\end{cases}
\end{equation*} Taking into account the following asymptotics: $g_H(\xi)\sim \sqrt{\frac{\pi}{2}}\xi^{H-1/2}e^{-\xi}$, $R_H(\xi)\sim \log(\xi)$  and $r_H(\xi)=1+o(1)$, as $\xi\rightarrow\infty$ (see Lemma \ref{LEM0-App}) we can write

\begin{equation*}
\begin{cases}
& \mathbf{P}_{1,1}  = \kappa  (-8y- 8x + 10u_0 )+o(y) \\
&\mathbf{P}_{2,1} = \kappa (56y- 112x + 70u_0 )+o(y) \\
&\mathbf{P}_{3,1} = \kappa (-16y- 28x + 28u_0 )+o(y) \\[3mm]
&\mathbf{P}_{1,2}  = 10u_0\kappa\sigma^2 \bm{\xi}  + 2\kappa\sigma^2 \left(  - 4Ax + 5u^* \right)+o(Ax) \\
&\mathbf{P}_{2,2} = 70u_0\kappa\sigma^2 \bm{\xi}   + 2\kappa\sigma^2 \left( - 56Ax + 35u^* \right)+o(Ax) \\
&\mathbf{P}_{3,2} = 28u_0\kappa\sigma^2 \bm{\xi}   + 2\kappa\sigma^2 \left( - 14Ax + 14u^* \right)+o(Ax) \\[3mm]
&\mathbf{P}_{1,3}  = \kappa\sigma^2 \bm{\zeta} (-8y- 8x + 10u_0 ) + 8\kappa\sigma^2\alpha \left( 2by + ax \right) +o(y)\\
&\mathbf{P}_{2,3} = \kappa\sigma^2 \bm{\zeta}(56y- 112x + 70u_0 ) - 8\kappa\sigma^2\alpha \left( 14by - 14ax \right)+o(y) \\ 
&\mathbf{P}_{3,3} = \kappa\sigma^2 \bm{\zeta}(-16y- 28x + 28u_0 ) - 2\kappa\sigma^2\alpha \left( - 16by - 14ax \right)+o(y)
\end{cases}
\end{equation*}With these arguments we evaluate the determinant of $\mathbf{P}$ and find that the first leading terms, such as, those of order $\OO(\alpha^l A^k x^j)$, $k,l=0,1$, $j=0,1,2,3$ are vanishing. So, we arrive at
\begin{eqnarray*} 
&& \abs{\mathbf{P}} \sim \kappa^3\sigma^4
\begin{vmatrix}
-8y- 8x + 10u_0    & 10(u_0\bm{\xi}+u^*) - 8Ax    &  ( 16b\alpha -8\bm{\zeta})y +(8a\alpha- 8\bm{\zeta})x +10u_0\bm{\zeta} \\
56y- 112x + 70u_0  & 70u_0(\bm{\xi}+u^*) - 112Ax  &   (-112b\alpha +56 \bm{\zeta} )y+(112a\alpha-112 \bm{\zeta} )x +70u_0\bm{\zeta}    \\
-16y- 28x + 28u_0  & 28(u_0\bm{\xi}+u^*)- 28Ax    & ( 32b\alpha -16\bm{\zeta})y+(28a\alpha - 28\bm{\zeta}  )x+ 28u_0\bm{\zeta} 
\end{vmatrix}\\
&& \sim -5376 u_0\kappa^3 \sigma^4 \alpha A xy,~~~(~A=R_H(\lambda\alpha)\sim\log(\lambda\alpha ),~x=g_H(\lambda\alpha),~y=g_H(2\lambda\alpha)~),
\end{eqnarray*}as $\alpha\rightarrow\infty$. Note that $u_0,\kappa,\sigma^2>0$. This means that we can find some $\widetilde{\alpha}>0$ (depending only on the parameter space $\Theta$), so that $\abs{\B \nabla \varrho_{\theta}(\alpha)}<0$, for all $\alpha>\widetilde{\alpha}$ (to specify $\widetilde{\alpha}$, just apply the same procedures as in Subsection \textbf{B.1}). We conclude that $\B \nabla \varrho_{\theta}(\alpha)$ is invertible, and by local univalent theorem the mapping $\theta \mapsto \B\varrho_{\theta}(\alpha)$ is locally injective. Therefore, the proof is complete.
\end{proof}

%
%
%
%
%
%
\subsection*{Disclosure statement}

Both authors have no conflicts of interest.
%
%
%
%
%
%
%
%
%
%
%
%
%
%

\subsection*{Acknowledgments}
The authors would like to express their gratitude to the anonymous reviewers for their valuable comments and constructive suggestions, which have significantly improved the quality of the paper.

\bibliographystyle{plain}
\bibliography{Biblioforall2}

\end{document}